\documentclass{amsart}

\usepackage{amsmath,amssymb,amsthm}
\usepackage{tikz-cd}
\usepackage[margin=1.15in]{geometry}

\newtheorem{theorem}{Theorem}[section]
\newtheorem{proposition}[theorem]{Proposition}
\newtheorem{lemma}[theorem]{Lemma}
\newtheorem{corollary}[theorem]{Corollary}
\newtheorem*{thmA}{Theorem A}
\newtheorem*{thmB}{Theorem B}
\newtheorem*{thmC}{Theorem C}
\theoremstyle{definition}
\newtheorem{definition}[theorem]{Definition}
\theoremstyle{remark}
\newtheorem{remark}[theorem]{Remark}

\newcommand{\A}{\mathbb{A}}
\newcommand{\RR}{\mathbb{R}}
\newcommand{\CC}{\mathbb{C}}
\newcommand{\ZZ}{\mathbb{Z}}
\newcommand{\cH}{\mathcal{H}}
\newcommand{\cK}{\mathcal{K}}
\newcommand{\RS}{\mathrm{RS}}
\newcommand{\co}{\mathrm{co}}
\newcommand{\parW}{\sigma_W}
\newcommand{\Atail}{\mathsf{B}}
\newcommand{\Apos}{\A_{+}}
\newcommand{\hp}{h_\partial}
\newcommand{\Pspl}{\widetilde\Pi}
\newcommand{\Ksm}{\cK_t}
\newcommand{\eigcut}{\lambda_0}
\newcommand{\Mf}{f}
\newcommand{\Mgen}[2]{x^{#1}_{#2}}
\newcommand{\Mhat}[1]{\hat x^{#1}}
\DeclareMathOperator{\sdet}{sdet}
\DeclareMathOperator{\Str}{Str}
\DeclareMathOperator{\Tr}{Tr}
\DeclareMathOperator{\rk}{rk}
\DeclareMathOperator{\sgn}{sgn}
\DeclareMathOperator{\End}{End}
\DeclareMathOperator{\Hom}{Hom}
\DeclareMathOperator{\im}{im}

\title{A Cheeger--M\"uller theorem for flat superconnections}
\author{Ryan Mickler}
\date{October 2026}
\dedicatory{For my teacher, Mathai Varghese}

\begin{document}

\begin{abstract}
For a flat superconnection ${\mathbb{A}}$ on a ${\mathbb{Z}}_2$-graded vector bundle $E$ over a closed oriented odd-dimensional manifold $M$, we prove a Cheeger-M\"uller type theorem as proposed by Arias Abad and Sch\"atz. Namely, that the topological Reidemeister norm of ${\mathbb{A}}$ as defined by them equals the Ray-Singer metric defined by Mathai and Wu. The main step proves an analytic reduction result identified by Arias Abad and Sch\"atz and formulated as a conjecture by the author, that the Ray-Singer metric coincides with that of the fibrewise cohomology $H(E,\partial)$ with its Gauss-Manin connection, where $\partial$ is the component of ${\mathbb{A}}$ of form degree zero, under the Knudsen-Mumford isomorphism associated to the form-degree spectral sequence. As an application, we consider Goette's Morse flat superconnection of a fibrewise Morse function on a fibre bundle with even-dimensional fibres over an oriented odd-dimensional base.
\end{abstract}
\maketitle

\section{Introduction}

Let $M$ be a closed oriented Riemannian manifold of odd dimension $n$, let $E=E^0\oplus E^1$ be a $\ZZ_2$-graded real vector bundle of finite rank over $M$ with a Euclidean metric, and let
\[
\A=\A_{[0]}+\A_{[1]}+\A_{[2]}+\cdots
\]
be a flat superconnection on $E$ in the sense of Quillen \cite{Q}, i.e. $\A$ is an odd operator on $\Omega(M,E)$ satisfying the Leibniz rule, $\A_{[j]}$ raises form degree by $j$, and $\A^2=0$. Thus $\A_{[1]}$ is a connection on $E$ that preserves the grading. For $j\ne1$, $\A_{[j]}$ is a $j$-form with values in the endomorphisms of $E$ of parity $j+1$ modulo $2$. Mathai and Wu \cite{MW2} defined the Ray--Singer metric $\|\cdot\|_{\RS}$ on the determinant line of $H(M,E,\A)$, the cohomology of $(\Omega(M,E),\A)$, and showed that it does not depend on the metrics when $n$ is odd.

Arias Abad and Sch\"atz \cite[Def.~3.8]{AS} constructed a Reidemeister norm $\tau_R$ on $\det H(M,E,\A)$ from a smooth triangulation, using the integration of flat superconnections to representations up to homotopy \cite{AS2,BS}. They proved its independence of the triangulation \cite[Thm.~3.9]{AS} and proposed its equality with the Ray--Singer metric \cite[\S1]{AS}. This would generalize the Cheeger--M\"uller theorem \cite{Ch,Mu,BZ}, the equality of analytic and Reidemeister torsion that was conjectured by Ray and Singer \cite{RS}. Our main result proves this equality. We work over $\RR$, as in \cite{AS} and \cite{Mi}. Complex coefficients are treated in Remark~\ref{rem:complex}.

\begin{thmA}[Cheeger--M\"uller theorem for flat superconnections]
Let $K$ be a smooth triangulation of $M$ with a total ordering of its vertices, let $M^K$ be the simplicial set of $K$, and let $\int E$ be the representation up to homotopy of $M^K$ obtained by integrating $\A$ \cite{AS2,BS}. Under the de Rham isomorphism $H(M,E,\A)\cong H(M^K,\textstyle\int E)$ of \cite[Cor.~2.14]{AS}, the Reidemeister norm $\tau_R$ of \cite{AS} equals the Ray--Singer metric of $(E,\A)$.
\end{thmA}

Arias Abad and Sch\"atz proved that their Reidemeister norm reduces to that of the fibrewise cohomology with its induced flat connection \cite[Cor.~3.12]{AS}. They observed that their proposed Cheeger--M\"uller theorem is therefore equivalent to the same reduction for the Mathai--Wu metric \cite[\S1]{AS}. Theorem~B below proves this analytic reduction, formulated as Conjecture~6.2 in \cite{Mi}, and hence gives Theorem~A.
Write $\partial=\A_{[0]}$, an odd bundle endomorphism with $\partial^2=0$. The constant-rank assumption in that conjecture is automatic: flatness makes $\partial$ parallel for the connection $\A_{[1]}$ (Lemma~\ref{lem:rank}). Hence $H_\partial=H(E,\partial)$ is a $\ZZ_2$-graded vector bundle, and $\A_{[1]}$ induces on it a flat connection $\nabla^{H}$, the Gauss--Manin connection \cite[Lemma~3.3]{AS}. Filtering $\Omega(M,E)$ by form degree gives a spectral sequence that converges to $H(M,E,\A)$. We write $E_r$ for its pages, so that $E_2=H(M,H_\partial,\nabla^{H})$. As in \cite[Rem.~3.10]{AS} and \cite[\S6]{Mi}, the Knudsen--Mumford isomorphisms of the pages, with the sign fixed in \S\ref{sec:conv}, compose to an isomorphism
\[
\kappa\colon \det H(M,H_\partial,\nabla^{H})\longrightarrow \det H(M,E,\A).
\]

\begin{thmB}
$\kappa$ is an isometry for the Ray--Singer metrics: $\|\kappa(\,\cdot\,)\|_{\RS(E,\A)}=\|\cdot\|_{\RS(H_\partial,\nabla^H)}$.
\end{thmB}

\subsection*{Strategy of the proof of Theorem~A}
Let $H(K;H_\partial)$ be the simplicial cohomology of $K$ with coefficients in the flat bundle $H_\partial$, and let $\kappa_K$ be defined in the same way as $\kappa$ from the spectral sequence of the cochains of $M^K$, filtered by cochain degree.

We consider the following square
\begin{equation}\label{eq:square}
\begin{tikzcd}[column sep=8.5em,row sep=4.5em]
\det H(M,H_\partial,\nabla^H)
  \arrow[r,"\kappa","\substack{\text{Theorem B}\\ \text{(conjectured in \cite{Mi})}}"']
  \arrow[d,"\mathrm{dR}","\substack{\text{Cheeger--M\"uller}\\ \text{\cite{BZ}, \cite[Thm.~6.2]{Fa}}}"']
& \det H(M,E,\A)
  \arrow[d,dashed,"\int"',"\text{Theorem A}"]
\\
\det H(K;H_\partial)
  \arrow[r,"\kappa_K","\substack{\text{Arias Abad--Sch\"atz}\\ \text{\cite[Cor.~3.12]{AS}}}"']
& \det H(M^K,\textstyle\int E)
\end{tikzcd}
\end{equation}

 The right vertical map $\int$ is the de Rham isomorphism of \cite[Cor.~2.14]{AS}, which comes from the $A_\infty$ de Rham theorem \cite{AS2,BS}. Define the left vertical map $\mathrm{dR}$ to be the determinant of its induced map on $E_2$. This agrees with the classical de Rham isomorphism for $H_\partial$ up to signs determined by the grading conventions (Proposition~\ref{prop:AS}). The map $\int$ is not a classical integration map. The differential on its target $C(M^K,\int E)$ has components of every cochain degree, given by $\partial$ at the vertices, parallel transport along the edges and higher homotopies on the higher simplices. The map $\int$ itself is built from iterated integrals of the components of $\A$ of positive form degree.
The top row has the Ray--Singer metrics, the lower left corner has Farber's Poincar\'e--Reidemeister metric of $H_\partial$ \cite{Fa}, and the lower right corner has $\tau_R$. Theorem~A says that the dashed arrow is an isometry. The left arrow is an isometry by the Cheeger--M\"uller theorem in the form of Bismut--Zhang and Farber \cite{BZ}, \cite[Thm.~6.2]{Fa}, since these signs do not change norms. The bottom arrow is an isometry by \cite[Cor.~3.12]{AS}, and the square commutes. Proposition~\ref{prop:AS} supplies the filtration comparisons needed for these two facts, including a direct proof of the compatibility with Poincar\'e duality.



\subsection*{Proof of Theorem~B}
The obstruction to extending the method of \cite{Mi} is the degree-zero component $\partial$. With $\rho_t=t^{N/2}$, where $N$ is the form-degree operator, the adiabatic method of \cite[\S2]{Mi} studies the rescaled family
\[
D_t=t^{-1/2}\rho_t\A\rho_t^{-1}=t^{-1/2}\partial+\A_{[1]}+t^{1/2}\A_{[2]}+t\A_{[3]}+\cdots,\qquad t>0.
\]
For $\A=d_H$ this is the deformation $d+\sum_{i\ge1}t^iH_{2i+1}$ of $d$, and perturbation theory at $t=0$ controls its small eigenvalues. For a general flat superconnection the term $t^{-1/2}\partial$ is singular and the family has no limit as $t\to0$. In \cite{Mi} only scalar fluxes are treated, for which $\partial=0$, and this term is the difficulty behind the closing remark of \cite[\S6]{Mi}, which asks for a more refined analysis of the eigenvalues of this family.

We avoid the singular family. The key idea is Hodge theory in the fibres of $E$. The metric on $E$ gives the decomposition $E=\cH\oplus E_c$, where $\cH=\ker\partial\cap\ker\partial^*\cong H_\partial$ is the bundle of $\partial$-harmonic elements and $E_c=\im\partial\oplus\im\partial^*$. In \S\ref{sec:split} we show that every flat superconnection is conjugate, by a gauge transformation $g$ built from the Green operator of $\partial$, to a direct sum $\A_\cH\oplus\A_c$. The harmonic part $\A_\cH$ is a flat superconnection on $\cH$ with no degree-zero part, whose connection part is the Gauss--Manin connection. The contractible part $\A_c$ is a flat superconnection on $E_c$ whose degree-zero part is $\partial|_{E_c}$, so that $(E_c,\partial)$ is exact in every fibre.

The same Hodge decomposition underlies the reduction of a regular representation up to homotopy to its cohomology \cite[Thm.~4.13]{AC}, and the splitting follows from the decomposition theorem for twisted objects \cite[Cor.~3.17]{CL} (Remark~\ref{rem:split}). On the combinatorial side, $\tau_R$ is invariant under fibrewise quasi-isomorphisms \cite[Def.~3.13, Prop.~3.15]{AS}. The conjugation preserves the Ray--Singer metric by a theorem of Mathai and Wu \cite[Thm.~6.1]{MW2}. Since $g$ has terms of positive form degree it need not preserve product metrics, and we reprove this by metric independence for arbitrary metrics on $\Lambda T^*M\otimes E$, which also gives the invariance of $\chi'$ (Corollary~\ref{cor:gauge}). The new point for torsion is the removal of $\A_c$. The Ray--Singer metric of $\A_c$ is trivial by a symmetry built from the Hodge star of $M$ (Proposition~\ref{prop:double}), and the spectral sequence of $\A_c$ vanishes from the first page on, so it does not affect $\kappa$. Thus all of $\partial$, and with it the singular term, lies in a summand that contributes nothing to either side of Theorem~B.

In \S\ref{sec:zero} we extend the method of \cite{Mi} to the harmonic part, comparing the small eigenvalues of $D_t$ with a finite model given by Hodge theory for $\nabla^H$ and homological perturbation. Comparing the scaling exponents also computes the derived Euler characteristic $\chi'(\A)$ (Definition~\ref{def:RS}), extending \cite[Thm.~2.9]{Mi}.

\begin{thmC}
Let $b_k(H^j_\partial)=\dim H^k(M,H^j_\partial)$ and let $d_r$ be the differentials of the spectral sequence. Then
\[
\chi'(\A)=\sum_{k}(-1)^k k\bigl(b_k(H^0_\partial)-b_k(H^1_\partial)\bigr)+\sum_{r\ge2}\bigl(\rk d_r|_{E_r^{0}}-\rk d_r|_{E_r^{1}}\bigr),
\]
where $E_r^{0}$ and $E_r^{1}$ are the even and odd parts of the $r$-th page.
\end{thmC}

For $M$ connected with $n\ge3$ and $H$ a non-zero multiple of the volume form, Theorem~C gives $\chi'(d+H)=\chi'(d)+1$, recovering \cite[Thm.~4.3]{Mi}.

As an application, let $p\colon X\to B$ be a fibre bundle with closed even-dimensional fibres over a closed oriented odd-dimensional manifold, with a flat vector bundle $F\to X$ and a fibrewise Morse function. Goette \cite{Go1,Go2} constructs from these data a flat superconnection $\A'$ on a finite-rank bundle $V\to B$ and a filtered quasi-isomorphism $I\colon(\Omega(X,F),d^F)\to(\Omega(B,V),\A')$. Theorem~B and Ma's functoriality theorem \cite{Ma} show that $I_*$ is an isometry from the Ray--Singer metric of $(X,F)$ to that of $(V,\A')$ (Corollary~\ref{cor:morse}). In \S\ref{ssec:doublecover} we work out an example over $S^3$ with fibre $S^1\times S^3$, in which the Morse differential and a degree-three term are both non-zero and the total space has $\ZZ/2$ in degree-four cohomology.

\subsection*{Relation to earlier work}
In \cite[Thm.~3.7]{Mi} the author proved Theorem~B for the scalar-flux complex $d+H$, with $H$ a closed sum of odd-degree forms of degrees at least three. That method combines Farber's small-eigenvalue analysis \cite{Fa95} with Forman's rescaling of filtered complexes \cite{Fo}. For a top-degree scalar flux on a unitary flat bundle, the torsion formula is due to Mathai and Wu \cite[Prop.~5.1]{MW1}. Sorto studied twisted torsion in the equivariant setting of finite group actions \cite{So}.

\section{Preliminaries}\label{sec:conv}

Vector bundles, vector spaces and determinant lines are real, except in Remark~\ref{rem:complex} and where a complexification is taken explicitly. A metric on a vector bundle means a Euclidean metric, and a norm on a line satisfies $\|zv\|=|z|\,\|v\|$.

\subsection{Determinant lines}
For a finite-dimensional $\ZZ_2$-graded vector space $V=V^0\oplus V^1$ put $\det V=\det V^0\otimes(\det V^1)^{-1}$. For an ordered basis $v=(v_1,\dots,v_m)$ of a vector space we write $\wedge v=v_1\wedge\dots\wedge v_m$ for the corresponding element of its top exterior power. A finite-dimensional $\ZZ_2$-graded complex $(C,\delta)$, with $\delta$ odd and $\delta^2=0$, has a Knudsen--Mumford isomorphism \cite{KM}
\[
\tau_{C,\delta}\colon \det C\longrightarrow \det H(C,\delta).
\]
As in \cite[App.~A, Rem.~3.10]{AS} and \cite[\S\S3, 6]{Mi}, these isomorphisms define the maps $\kappa$ below, and we fix their sign as follows. For $i=0,1$ choose ordered sets $b^i\subset C^i$ such that $\delta b^i$ is a basis of $\delta(C^i)\subset C^{i+1}$, and $z^i\subset\ker\delta\cap C^i$ whose classes $[z^i]$ form a basis of $H^i(C,\delta)$. Then $(\delta b^1,z^0,b^0)$ is a basis of $C^0$ and $(\delta b^0,z^1,b^1)$ is a basis of $C^1$, and we put
\[
\tau_{C,\delta}\bigl(\wedge(\delta b^1,z^0,b^0)\otimes\wedge(\delta b^0,z^1,b^1)^{-1}\bigr)=\wedge[z^0]\otimes\wedge[z^1]^{-1}.
\]
This does not depend on the choices. Changing $b^i$ leaves the argument of $\tau_{C,\delta}$ unchanged, and changing $z^i$ multiplies both sides by the same determinant. For $C^0=\RR e_0$ and $C^1=\RR\,\delta e_0$ with $\delta e_0\neq0$ it gives
\[
\tau_{C,\delta}\bigl(e_0\otimes(\delta e_0)^{-1}\bigr)=1\in\RR=\det H(C,\delta).
\]
Property~(KM2) below and Lemma~\ref{lem:finiteRS} depend on this normalization. If $V$ has a finite filtration with associated graded space $\operatorname{gr}V$, lifting an ordered basis of $\operatorname{gr}V$, adapted to its grading, to $V$ gives an isomorphism $\det\operatorname{gr}V\to\det V$. It does not depend on the lifts. We use three properties of $\tau$.
\begin{enumerate}
\item[(KM1)] If $\phi\colon(C,\delta)\to(C',\delta')$ is an isomorphism of complexes, then $\tau_{C',\delta'}\circ\det\phi=H(\phi)_*\circ\tau_{C,\delta}$.
\item[(KM2)] For $c>0$, $\tau_{C,c\delta}=c^{\,r(\delta)}\,\tau_{C,\delta}$, where $r(\delta)=\rk(\delta|_{C^0})-\rk(\delta|_{C^1})$.
\item[(KM3)] Let $C$ be filtered, with $\delta$ preserving the filtration, and let $E_r$ and $d_r$ be the pages and differentials of its spectral sequence. Put $n_r=\rk(d_r|_{E_r^0})$, $m_r=\rk(d_r|_{E_r^1})$ and $\sgn(C,\delta)=(-1)^{\sum_{r<s}(n_rn_s+m_rm_s)}$. Then $\sgn(C,\delta)\,\tau_{C,\delta}$ is the composite of $\det C=\det E_0$ with the isomorphisms $\tau_{E_r,d_r}\colon\det E_r\to\det E_{r+1}$ and $\det E_\infty=\det\operatorname{gr}H(C,\delta)\to\det H(C,\delta)$. Compare \cite[Thm.~1.3]{Fd} and \cite[Prop.~A.7]{AS}, where this is stated up to sign.
\end{enumerate}
(KM1) and (KM2) follow from the definition of $\tau$. For (KM3), choose $b^i_r\subset C^i$ representing a basis, homogeneous for the filtration degree, of a complement of $\ker d_r$ in $E_r^i$, and $z^i$ representing a homogeneous basis of $E_\infty^i$. Then $\delta b^i_r$ represents a basis of $d_r(E_r^i)$, and by descending induction on $r$ the vectors of index at least $r$ represent a basis of $E_r$. The composite in (KM3) is therefore computed by the bases $(\delta b^1_0,\delta b^1_1,\dots,z^0,\dots,b^0_1,b^0_0)$ of $C^0$ and $(\delta b^0_0,\delta b^0_1,\dots,z^1,\dots,b^1_1,b^1_0)$ of $C^1$. The definition of $\tau_{C,\delta}$ requires $\delta b^i$ in the same order as $b^i$. This reverses the blocks $\delta b^i_r$ and gives the sign $\sgn(C,\delta)$.

For a filtered complex we write $\kappa\colon\det E_2\to\det H(C,\delta)$ for the composite of the isomorphisms $\tau_{E_r,d_r}$ with $r\ge2$ and $\det E_\infty\to\det H(C,\delta)$. This is the map of \cite[Rem.~3.10]{AS} and \cite[\S\S3, 6]{Mi}, with our sign convention. A filtered chain map that induces isomorphisms of the pages $E_r$ for $r\ge2$ intertwines the two maps $\kappa$ exactly, by (KM1) on each page. The sign in (KM3) enters only through Lemma~\ref{lem:model}.

\begin{lemma}\label{lem:finiteRS}
Let $(C,\delta)$ be a finite-dimensional $\ZZ_2$-graded complex with a metric for which $C^0\perp C^1$. Put $\sdet'(\delta^*\delta)=\det'(\delta^*\delta|_{C^0})/\det'(\delta^*\delta|_{C^1})$, where $\det'$ is the product of the non-zero eigenvalues, and let $|\cdot|_{\mathrm{harm}}$ be the metric on $\det H(C,\delta)$ induced by harmonic representatives. Then
\[
\|\cdot\|_{C,\delta}=\sdet'(\delta^*\delta)^{-1/2}\,|\cdot|_{\mathrm{harm}}
\]
is the metric for which $\tau_{C,\delta}$ is an isometry from $\det C$, with the metric induced by the metric on $C$, to $\det H(C,\delta)$.
\end{lemma}

\begin{proof}
Write $C=\cK\oplus\im\delta\oplus\im\delta^*$ with $\cK=\ker\delta\cap\ker\delta^*$. On $\im\delta^*\cap C^i$ the map $\delta$ is injective with singular values the square roots of the non-zero eigenvalues of $\delta^*\delta|_{C^i}$. Computing $\tau_{C,\delta}$ in an orthonormal basis adapted to this decomposition gives the claim. For $C^0=\RR\xrightarrow{\lambda}\RR=C^1$ with the standard metrics, both descriptions give $\|1\|=|\lambda|^{-1}$ by the normalization of $\tau$.
\end{proof}

\subsection{The Ray--Singer metric}
Let $\Omega(M,E)^i=\bigoplus_{k+j\equiv i}\Omega^k(M,E^j)$, so that $\A$ is odd. Fix a metric on $\Lambda T^*M\otimes E$ for which the even and odd parts are orthogonal. Such a metric need not come from metrics on $TM$ and $E$. We use this generality in Corollary~\ref{cor:gauge}. Let $\A^*$ be the formal adjoint for the $L^2$ inner product defined by this metric and the Riemannian density, and let $\Delta=\A\A^*+\A^*\A$. The principal symbol of $\A$ is that of $\A_{[1]}$, which is exterior multiplication by $\xi\in T^*M$, up to a factor $\mathrm{i}$. For $0\neq\xi\in T^*_xM$ the complex $(\Lambda T^*_xM\otimes E_x,\xi\wedge)$ is exact, so the principal symbol of $\Delta$ is positive definite for every choice of the metric, and $\Delta$ is elliptic. Let $P$ be the orthogonal projection onto $\ker\Delta\cong H(M,E,\A)$ and $\Pi^{\co}$ that onto the closure of $\im\A^*$.

\begin{definition}\label{def:RS}
For $i=0,1$ and $\operatorname{Re}s\gg0$ let $\zeta_i(s)=\Tr(\Pi^{\co}\Delta^{-s}|_{\Omega(M,E)^i})$. These are the zeta functions of $\A^*\A$ on even and odd forms, with the kernel removed. They are regular at $s=0$ (\cite[Prop.~2.1]{MW2} and Lemma~\ref{lem:regular}). The \emph{Ray--Singer metric} on $\det H(M,E,\A)$ is
\[
\|\cdot\|_{\RS}=\exp\bigl(\tfrac12\zeta_0'(0)-\tfrac12\zeta_1'(0)\bigr)\,|\cdot|_{L^2},
\]
where $|\cdot|_{L^2}$ is induced by the harmonic representatives. We also define the \emph{derived Euler characteristic}
\[
\chi'(\A)=\zeta_0(0)-\zeta_1(0),
\]
as in the scalar-flux case of \cite[\S1.2]{Mi}. For an ordinary flat bundle, Theorem~C identifies it with $\sum_q(-1)^qq\dim H^q$, the quantity used in \cite{BZ}.
\end{definition}

With $\sdet'(\A^*\A)=\exp(-\zeta_0'(0)+\zeta_1'(0))$ the definition reads $\|\cdot\|_{\RS}=\sdet'(\A^*\A)^{-1/2}|\cdot|_{L^2}$.\footnote{The displayed definition of $\|\cdot\|_H$ in \cite[\S1.1]{Mi} has the exponent $+1/2$. On our determinant line it must be $-1/2$, as Lemma~\ref{lem:finiteRS} shows. The positive exponent is correct on the dual line $(\det H)^{-1}$.} This is the metric for which the torsion element of Mathai and Wu \cite[Def.~2.2]{MW2} has norm one, and by Lemma~\ref{lem:finiteRS} it is the metric $\|\cdot\|_{C,\delta}$ for a finite complex.

The same definition applies to any odd first-order differential operator $D$ on $\Omega(M,E)$ with $D^2=0$ and $\Delta_D=DD^*+D^*D$ elliptic, for instance $D=c\A$ with $c>0$, for which
\begin{equation}\label{eq:scale}
\|\cdot\|_{\RS(c\A)}=c^{-\chi'(\A)}\,\|\cdot\|_{\RS(\A)},\qquad \chi'(c\A)=\chi'(\A),
\end{equation}
since $\zeta_i$ for $c\A$ equals $c^{-2s}\zeta_i(s)$ and the harmonic forms are unchanged.

\begin{lemma}\label{lem:regular}
Let $D$ be as above and $n$ odd.
\begin{enumerate}
\item $\zeta_i(s)=\Gamma(s+1)^{-1}\int_0^\infty v^{s}\,\Tr\bigl(D^*D\,e^{-v\Delta}|_{\Omega(M,E)^i}\bigr)\,dv$, and $\zeta_i$ is regular at $s=0$.
\item For every bundle endomorphism $\theta$ of $\Lambda T^*M\otimes E$ preserving the parity, $\Str(\theta(1-P)\Delta^{-s})$ is regular at $s=0$ with value $-\Str(\theta P)$.
\end{enumerate}
\end{lemma}

\begin{proof}
Since $\Pi^{\co}\Delta=D^*D$, we have $\tfrac{d}{dv}\Tr(\Pi^{\co}e^{-v\Delta}|_{\Omega(M,E)^i})=-\Tr(D^*De^{-v\Delta}|_{\Omega(M,E)^i})$, and inserting this into the Mellin transform gives the formula in (1). For a differential operator $Q$ of order $m$ the heat trace has an expansion $\Tr(Qe^{-v\Delta})\sim\sum_{j\ge0}a_j(Q)\,v^{(j-n-m)/2}$ as $v\to0$, with $a_j(Q)=0$ for $j+m$ odd, because the homogeneous parts of the symbols of $Q$ and of the parametrix of $\Delta$ are even or odd in $\xi$ according to their degree \cite[Lemma~1.7.7]{Gi}. A pole of $\zeta_i$ at $s=0$ would come from $a_n(D^*D)$, and $n+2$ is odd. For (2), write
\[
\Str\bigl(\theta(1-P)\Delta^{-s}\bigr)=\Gamma(s)^{-1}\int_0^\infty v^{s-1}\,\Str\bigl(\theta(e^{-v\Delta}-P)\bigr)\,dv.
\]
The integral over $(1,\infty)$ is entire in $s$, and it is killed by the zero of $\Gamma(s)^{-1}$ at $s=0$. In the integral over $(0,1)$, the coefficient of $v^0$ in $\Str(\theta e^{-v\Delta})$ is $a_n(\theta)=0$ because $n$ is odd, and the constant term $-\Str(\theta P)$ gives the pole $-\Str(\theta P)/s$. Hence the value at $s=0$ is $-\Str(\theta P)$.
\end{proof}

Several results we cite, notably \cite{BZ,Ma}, are stated for complex flat bundles. We apply them to real bundles through complexification. We equip $E_\CC=E\otimes\CC$ with the complex-linear extension $\A_\CC$ of $\A$ and the Hermitian extension of the metric.

\begin{lemma}\label{lem:cplx}
$H(M,E_\CC,\A_\CC)=H(M,E,\A)\otimes\CC$, and under $\det_\CC H(M,E_\CC,\A_\CC)=\det H(M,E,\A)\otimes\CC$ the Ray--Singer metric of $(E_\CC,\A_\CC)$ restricts to that of $(E,\A)$, and $\chi'(\A_\CC)=\chi'(\A)$. Likewise, complexification preserves the maps $\tau_{C,\delta}$ and $\kappa$, the metrics $\|\cdot\|_{C,\delta}$, and Farber's Poincar\'e--Reidemeister metric of a flat bundle.
\end{lemma}

\begin{proof}
The Laplacian of $\A_\CC$ is the complex-linear extension of that of $\A$. It has the same eigenvalues with the same multiplicities, so the zeta functions agree, and a real orthonormal basis of harmonic forms is a complex orthonormal basis. For a finite real complex, real choices of $b^i$ and $z^i$ in the definition of $\tau_{C,\delta}$ are also valid complex choices. Farber's metric is assembled from such finite complexes and the Poincar\'e pairing, both of which are compatible with complexification.
\end{proof}

\subsection{Metric independence}\label{sec:metric}

We use metric independence in the general elliptic-complex setting of Mathai and Wu \cite[Thm.~4.1]{MW2}, with coefficient bundle $\Lambda T^*M\otimes E$. Allowing arbitrary parity-preserving metrics on this bundle is important: a gauge transformation with positive-form-degree terms need not preserve product metrics. We recall the variation argument, including the invariance of $\chi'$.

\begin{proposition}\label{prop:metric}
Let $n$ be odd. The Ray--Singer metric and $\chi'(\A)$ do not depend on the choice of a metric on $\Lambda T^*M\otimes E$ for which the even and odd parts are orthogonal.
\end{proposition}

\begin{proof}
Metrics on $\Lambda T^*M\otimes E$ for which the even and odd parts are orthogonal form a convex set, so any two of them are joined by a smooth path $\langle\cdot,\cdot\rangle_u$. The $L^2$ inner product is formed with a fixed smooth density on $M$, and a change of density can be absorbed into the fibre metric. The proposition therefore includes independence of the Riemannian metric. Write $\langle a,b\rangle_u=\langle g_ua,b\rangle_0$, with $g_u$ an even positive bundle endomorphism, and put $\theta=g_u^{-1}\tfrac{d}{du}g_u$. Then $\tfrac{d}{du}(\A^*)=-\theta\A^*+\A^*\theta$. The dimension of $\ker\Delta\cong H(M,E,\A)$ does not depend on $u$, so $P$ depends smoothly on $u$. Perturbing the non-zero eigenvalues of $\A^*\A$ gives
\[
\frac{\partial}{\partial u}\bigl(\zeta_0(s)-\zeta_1(s)\bigr)=s\,\Str\bigl(\theta(1-P)\Delta^{-s}\bigr).
\]
The two terms of $\tfrac{d}{du}(\A^*)$ contribute $s\,\Str(\theta\Pi^{\co}\Delta^{-s})$ and $s\,\Str(\theta\,\Pi^{\mathrm{ex}}\Delta^{-s})$, where $\Pi^{\mathrm{ex}}=1-P-\Pi^{\co}$ is the projection onto the closure of $\im\A$. By Lemma~\ref{lem:regular}(2) the right side vanishes at $s=0$, so $\tfrac{d}{du}\chi'(\A)=0$, and its derivative $\tfrac{\partial}{\partial s}$ at $s=0$ is $-\Str(\theta P)$. Since the harmonic representative of a class minimizes the norm in its class, $\tfrac{d}{du}\log|\cdot|_{L^2}=\tfrac12\Str(\theta P)$. Hence $\tfrac{d}{du}\log\|\cdot\|_{\RS}=-\tfrac12\Str(\theta P)+\tfrac12\Str(\theta P)=0$.
\end{proof}

\begin{corollary}\label{cor:gauge}
Let $g\in\Omega(M,\End E)$ be even and invertible. Then $g$ induces an isometry
\[
\det H(M,E,\A)\longrightarrow\det H(M,E,g\A g^{-1})
\]
for the Ray--Singer metrics, and $\chi'(g\A g^{-1})=\chi'(\A)$.
\end{corollary}

\begin{proof}
$g$ is an isometry from $(\Omega(M,E),\A)$ with the metric $\langle g\cdot,g\cdot\rangle$ to $(\Omega(M,E),g\A g^{-1})$ with the original metric. Apply Proposition~\ref{prop:metric}.
\end{proof}

Mathai and Wu prove the isometry for $g$ joined to $1$ by a path of even invertible elements \cite[Thm.~6.1]{MW2}. This covers $g=1+\eta$ by the path $1+u\eta$, $u\in[0,1]$.

\section{The Hodge splitting}\label{sec:split}

Using Hodge theory for $\partial$ in the fibres of $E$, we split $\A$ by a gauge transformation into a harmonic part, which is a flat superconnection on $H_\partial$ with no degree-zero part, and a contractible part. We then show that the contractible part is invisible to both sides of Theorem~B.

\begin{lemma}\label{lem:rank}
$\partial$ has constant rank on each component of $M$.
\end{lemma}

\begin{proof}
The component of $\A^2=0$ of form degree one is $[\A_{[1]},\partial]=0$, one of the relations in \cite[Rem.~2.2]{AS}. Thus $\partial$ is parallel for the connection $\A_{[1]}$.
\end{proof}

Using the metric on $E$, let $\Delta_\partial=\partial\partial^*+\partial^*\partial$ be the fibrewise Laplacian. By Lemma~\ref{lem:rank} it has constant rank, so the fibrewise Hodge decomposition $E=\cH\oplus E_c$, with $\cH=\ker\partial\cap\ker\partial^*$ and $E_c=\im\partial\oplus\im\partial^*$, is smooth. Moreover $\cH\cong H_\partial$, and $(E_c,\partial)$ is acyclic. Let $\Pi_\cH$ be the projection onto $\cH$, let $G_\partial$ be the Green operator, which is zero on $\cH$ and inverts $\Delta_\partial$ on $E_c$, and put $\hp=G_\partial\partial^*$. Odd endomorphisms $T$ of $E$, in particular $\partial$ and $\hp$, act on $\Omega(M,E)$ with the Koszul sign, $T(\alpha\otimes e)=(-1)^{k}\alpha\otimes Te$ for a $k$-form $\alpha$. Then
\begin{equation}\label{eq:fibrehodge}
1-\Pi_\cH=\partial\hp+\hp\partial,\qquad \hp^2=0,\qquad \hp\Pi_\cH=\Pi_\cH\hp=0.
\end{equation}
Under $\cH\cong H_\partial$ the Gauss--Manin connection is $\nabla^H=\Pi_\cH\A_{[1]}|_\cH$.

\begin{proposition}\label{prop:split}
Put $\Apos=\A-\partial$, and let
\[
\Pspl=(1+\hp\Apos)^{-1}\,\Pi_\cH\,(1+\Apos\hp)^{-1},
\qquad
g=\Pi_\cH\Pspl+(1-\Pi_\cH)(1-\Pspl).
\]
Then $g=1+\eta$ with $\eta\in\bigoplus_{k\ge1}\Omega^k(M,\End E)$ even, and $g\A g^{-1}=\A_\cH\oplus\A_c$ with respect to $E=\cH\oplus E_c$. Here $\A_\cH$ is a flat superconnection on $\cH$ with $\A_{\cH,[0]}=0$ and $\A_{\cH,[1]}=\nabla^H$, given by $\A_\cH=\Pi_\cH\Apos(1+\hp\Apos)^{-1}$ on $\Omega(M,\cH)$, and $\A_c$ is a flat superconnection on $E_c$ with $\A_{c,[0]}=\partial|_{E_c}$.
\end{proposition}

\begin{proof}
Since $\hp$ and $\Pi_\cH$ preserve form degree and $\Apos$ raises it, $\hp\Apos$ and $\Apos\hp$ are nilpotent, so these are finite sums. Let $\iota\colon\Omega(M,\cH)\to\Omega(M,E)$ be the inclusion. By \eqref{eq:fibrehodge}, $(\iota,\Pi_\cH,\hp)$ is a special deformation retract of $(\Omega(M,E),\partial)$ onto $\Omega(M,\cH)$ with zero differential. Put $\iota_\infty=(1+\hp\Apos)^{-1}\iota$, $\pi_\infty=\Pi_\cH(1+\Apos\hp)^{-1}$, $k=\hp(1+\Apos\hp)^{-1}$ and $\A'=\Pi_\cH\Apos\iota_\infty$. The homological perturbation lemma \cite[Lemma~2.4 and Rem.~2.3]{Cr}, with $\hp$ replaced by $-\hp$ as in Lemma~\ref{lem:model}, gives $\pi_\infty\iota_\infty=1$, $\A\iota_\infty=\iota_\infty\A'$ and $\Pspl=\iota_\infty\pi_\infty=1-\A k-k\A$. Hence $\Pspl^2=\Pspl$ and $\Pspl$ commutes with $\A$.

We spell out the $\Omega(M)$-linearity, as in the proof of \cite[Thm.~4.13]{AC}. For a homogeneous scalar form $\alpha$, the Leibniz rule gives $[\Apos,\alpha]=d\alpha$, with the graded commutator. Commute $\alpha$ through the finite series defining $\iota_\infty$ and $\pi_\infty$. Every derivative term contains $\hp(d\alpha)\hp$, $\hp(d\alpha)\iota$ or $\Pi_\cH(d\alpha)\hp$. Moving $d\alpha$ past the adjacent fibrewise maps makes these terms vanish by $\hp^2=0$, $\hp\iota=0$ and $\Pi_\cH\hp=0$, respectively. Hence $\iota_\infty$ and $\pi_\infty$ are $\Omega(M)$-linear, and so is $\Pspl=\iota_\infty\pi_\infty$. Thus $\Pspl$ is an even element of $\Omega(M,\End E)$, and $\Pspl-\Pi_\cH$ has positive form degree. Therefore $g=1+\eta$ as stated, and $g$ is invertible because $\eta$ is nilpotent.

From $\Pi_\cH g=\Pi_\cH\Pspl=g\Pspl$ we get $g\Pspl g^{-1}=\Pi_\cH$. So $g\A g^{-1}$ commutes with $\Pi_\cH$, which gives the direct sum $\A_\cH\oplus\A_c$ of flat superconnections. Its component of form degree zero is $\partial$, and its component of form degree one is $\A_{[1]}+[\eta_{[1]},\partial]$, where $\eta_{[1]}$ is the component of $\eta$ of form degree one. Since $\partial\Pi_\cH=\Pi_\cH\partial=0$, this gives $\A_{\cH,[0]}=0$, $\A_{c,[0]}=\partial|_{E_c}$ and $\A_{\cH,[1]}=\nabla^H$. Finally, $\Pspl\iota_\infty=\iota_\infty$ and $\Pi_\cH\hp=0$ give $g\iota_\infty=\Pi_\cH\iota_\infty=\iota$. Hence $g\A g^{-1}\iota=g\A\iota_\infty=g\iota_\infty\A'=\iota\A'$, that is, $\A_\cH=\A'$.
\end{proof}

Write $E_c=W\oplus\partial W$ with $W=\im\partial^*$. Choose a metric on $W$ with $W^0\perp W^1$, transport it to $\partial W$ by $\partial$, and make $W\perp\partial W$. By Proposition~\ref{prop:metric} this does not change the Ray--Singer metric. Let $\parW$ be the parity operator of $\Omega(M,W)$, equal to $(-1)^{k+j}$ on $\Omega^k(M,W^j)$. We identify $\Omega(M,W)\oplus\Omega(M,W)$ with $\Omega(M,E_c)$ by $(x,y)\mapsto x+\partial \parW y$. For a metric connection $\nabla$ on $W$ that preserves the grading of $W$, with curvature $R=\nabla^2$, define
\[
D_\nabla\begin{pmatrix}x\\ y\end{pmatrix}=\begin{pmatrix}\nabla x-R\parW y\\ \parW x+\nabla y\end{pmatrix}.
\]
Using $\parW\nabla=-\nabla \parW$, $\parW R=R\parW$ and $\nabla R=R\nabla$ one checks $D_\nabla^2=0$. $D_\nabla$ is a flat superconnection on $E_c$: its degree-zero part is $\partial$, its degree-one part is $\nabla\oplus\nabla$, and its degree-two part is $-R$ composed with the inverse of $\partial$. This is the ``double'' of \cite[Ex.~3.8]{AC}.

\begin{proposition}\label{prop:normal}
$\A_c$ is gauge equivalent to $D_\nabla$ by an even $g=1+\eta$ with $\eta\in\bigoplus_{k\ge1}\Omega^k(M,\End E_c)$.
\end{proposition}

\begin{proof}
The complex $\End E_c$ with differential $[\partial,\cdot\,]$ is acyclic. If $g\A_cg^{-1}$ and $D_\nabla$ agree in form degrees less than $k$, the components of form degree $k$ of $(g\A_cg^{-1})^2=0$ and $D_\nabla^2=0$ show that their difference $X_k$ in degree $k$ satisfies $[\partial,X_k]=0$. Since $\partial$ has constant rank, $X_k=[\partial,Y_k]$ with a smooth even $Y_k\in\Omega^k(M,\End E_c)$. Replacing $g$ by $(1+Y_k)g$ leaves the components of $g\A_cg^{-1}$ of form degree less than $k$ unchanged and subtracts $X_k$ from the component of degree $k$. For $k=1$ this changes the connection on $W$ to $\nabla$. The induction terminates after degree $n$.
\end{proof}

\begin{remark}\label{rem:split}
The operator $\A_\cH$ is the structure on the cohomology bundle constructed by Arias Abad and Crainic \cite[Thm.~4.13]{AC}, from the same fibrewise Hodge decomposition and the homological perturbation lemma. Chuang and Lazarev give a conjugation of this kind for a complex with an abstract Hodge decomposition \cite[Cor.~3.13]{CL}\footnote{In \cite[Cor.~3.13]{CL}, with $t$ the projection, $s$ the homotopy, $x$ the perturbation, $\alpha=(1+sx)^{-1}$ and $\beta=(1+xs)^{-1}$, the term $t\alpha xt$ should read $t\beta xt$, as in \cite[Thm.~3.11]{CL}. Since $t\alpha=t$, the printed term is $txt$.}, and deduce the decomposition theorem for twisted objects \cite[Cor.~3.17]{CL}. Applied to $(\Omega(M,E),\partial)$ with the perturbation $\Apos$, their conjugation has the same idempotent $\Pspl$ and the same harmonic part, and it sends $\A$ to $\A_\cH\oplus\partial|_{E_c}$. That conjugation is not $\Omega(M)$-linear, since $\partial|_{E_c}$ does not satisfy the Leibniz rule.

An $\Omega(M)$-linear conjugation comes from the same results with the base $D_0=\nabla^H\oplus D_\nabla$ on $\cH\oplus E_c$. The identities
\[
D_0\hp+\hp D_0=1-\Pi_\cH,\qquad [D_0,\Pi_\cH]=0
\]
and \eqref{eq:fibrehodge} show that $\hp$ and $\Pi_\cH$ form an abstract Hodge decomposition for $D_0$. The perturbation $\A-D_0$ lies in $\bigoplus_{k\ge1}\Omega^k(M,\End E)$. Moreover, $B_0=D_0-\partial$ preserves the splitting and satisfies
\[
\hp B_0\hp=0,\qquad \hp B_0\Pi_\cH=\Pi_\cH B_0\hp=0.
\]
Thus replacing $\Apos$ by $\A-D_0$ in the finite perturbation series gives the same projector $\Pspl$ and, after adding the base differential $\nabla^H$, the same harmonic operator $\A_\cH$. Hence \cite[Thm.~3.11, Cor.~3.13]{CL}, whose proofs are computations with these relations and apply to $\ZZ_2$-gradings, give an even $g'=1+\eta'$ with $\eta'\in\bigoplus_{k\ge1}\Omega^k(M,\End E)$, $g'\A g'^{-1}=\A_\cH\oplus D_\nabla$ and $g'^{-1}\Pi_\cH g'=\Pspl$. This combines Propositions~\ref{prop:split} and~\ref{prop:normal} in one step.

Such splittings appear in \cite[Thm.~6.6, Prop.~6.3]{GSM} for two-term flat Lie algebroid superconnections, and in the proof of \cite[Prop.~6.7]{CHL} for reduced free twisted modules over a non-negatively graded pseudo-compact dg algebra. For a trivial finite-rank bundle with constant $\partial$, the latter argument applies to $\Omega(M)$ without completion, because positive form degree makes the perturbation series finite.
\end{remark}

\begin{proposition}\label{prop:double}
Let $n$ be odd. Then the complex $(\Omega(M,E_c),D_\nabla)$ is acyclic, its Ray--Singer metric is the standard metric on $\det 0=\RR$, and $\chi'(D_\nabla)=0$.
\end{proposition}

\begin{proof}
Acyclicity follows from the spectral sequence of the form-degree filtration, whose first page is $\Omega(M,H(E_c,\partial))=0$. Let $\star$ be the Hodge star defined by the metric and the orientation of $M$, and put $\tilde\star=(-1)^{k(k+1)/2}\star$ on $k$-forms, extended to $\Omega(M,W)$. Then $\tilde\star$ is orthogonal. We use that $n$ is odd three times. First, $\star$ maps $k$-forms to $(n-k)$-forms and so reverses the parity of the form degree, which gives $\tilde\star \parW=-\parW\tilde\star$. Second, $\star\star=(-1)^{k(n-k)}=1$, and the adjoint of the metric connection $\nabla$ on $k$-forms is $\nabla^*=(-1)^{n(k+1)+1}\star\nabla\star=(-1)^k\star\nabla\star$. Third, for a $k$-form $\omega$ these give
\[
\tilde\star\nabla\omega=(-1)^{(k+1)(k+2)/2}\star\nabla\omega,
\qquad
\nabla^*\tilde\star\omega=(-1)^{k(k+1)/2+n-k}\star\nabla\omega,
\]
and the two exponents differ by $2k+1-n$, which is even. Hence $\tilde\star\nabla=\nabla^*\tilde\star$. Squaring this identity gives the corresponding identity for the curvature:
\[
\tilde\star R=\tilde\star\nabla^2=(\nabla^*)^2\tilde\star=R^*\tilde\star.
\]
With
\[
J=\begin{pmatrix}0&-\tilde\star\\ \tilde\star&0\end{pmatrix}
\]
these identities and $\parW R^*=R^*\parW$ give $JD_\nabla=D_\nabla^*J$. $J$ is orthogonal and preserves the parity of $\Omega(M,E_c)$, so $D_\nabla^*D_\nabla$ on even forms is orthogonally equivalent to $D_\nabla D_\nabla^*$ on even forms. The latter has the same non-zero spectrum as $D_\nabla^*D_\nabla$ on odd forms. Therefore $\zeta_0=\zeta_1$.
\end{proof}

\begin{corollary}\label{cor:reduction}
Theorem~B holds for $(E,\A)$ if it holds for $(\cH,\A_\cH)$, and $\chi'(\A)=\chi'(\A_\cH)$.
\end{corollary}

\begin{proof}
The gauge transformations of Propositions~\ref{prop:split} and~\ref{prop:normal} preserve the form-degree filtration and induce the identity on the first page $\Omega(M,H_\partial)$, so they intertwine the maps $\kappa$. The spectral sequence of $\A_c$ vanishes from the first page on, so $\kappa_{\A_\cH\oplus\A_c}=\kappa_{\A_\cH}$. By Proposition~\ref{prop:metric} we may take $\cH\perp E_c$, and the Ray--Singer metric of an orthogonal direct sum is the tensor product. Corollary~\ref{cor:gauge} and Proposition~\ref{prop:double} give the claim.
\end{proof}

\section{The case $\partial=0$}\label{sec:zero}

We extend the results of \cite{Mi}, proved there for the trivial line bundle, to non-trivial bundles. Let $\cH$ be a $\ZZ_2$-graded vector bundle with a flat connection $\nabla$, not necessarily compatible with a metric, and let $\A=\nabla+\Atail$ be a flat superconnection with $\Atail=\sum_{j\ge2}\A_{[j]}$. Throughout this section we fix the product metric on $\Lambda T^*M\otimes\cH$ induced by a Riemannian metric on $M$ and a metric on $\cH$ with $\cH^0\perp\cH^1$. Proposition~\ref{prop:metric} permits this choice without changing the quantities to be compared. In particular, distinct form degrees are orthogonal, $\nabla^*$ lowers form degree by one, and $\Delta_\nabla$ preserves form degree.

\begin{theorem}\label{thm:zero}
Let $n$ be odd. Then $\kappa\colon\det H(M,\cH,\nabla)\to\det H(M,\cH,\A)$ is an isometry for the Ray--Singer metrics, and
\[
\chi'(\A)=\sum_k(-1)^kk\bigl(b_k(\cH^0)-b_k(\cH^1)\bigr)+\sum_{r\ge2}\bigl(\rk d_r|_{E_r^{0}}-\rk d_r|_{E_r^{1}}\bigr),
\]
where $b_k(\cH^j)=\dim H^k(M,\cH^j)$ and $d_r$ are the differentials of the spectral sequence.
\end{theorem}

\subsection{A finite model}
Let $\Delta_\nabla=\nabla\nabla^*+\nabla^*\nabla$, let $\cK=\ker\Delta_\nabla\cong H(M,\cH,\nabla)$, let $i\colon\cK\to\Omega(M,\cH)$ be the inclusion and $\Pi_\cK$ the orthogonal projection, and let $G$ be the Green operator. With $h=\nabla^*G$ we have $1-i\Pi_\cK=\nabla h+h\nabla$, $\Pi_\cK h=0$, $hi=0$ and $h^2=0$. The operator $G$ preserves form degree and $h$ lowers it by one. Thus $h\Atail$ raises form degree by at least one, and
\[
(1+h\Atail)^{-1}=\sum_{m=0}^{n}(-h\Atail)^m.
\]
All homological perturbation expressions below are therefore finite sums.

\begin{lemma}\label{lem:model}
Put $\delta=\Pi_\cK\Atail(1+h\Atail)^{-1}i$ and $i_\infty=(1+h\Atail)^{-1}i$. Then $\delta$ is an odd operator on $\cK$ with $\delta^2=0$ that raises form degree by at least two, $i_\infty\colon(\cK,\delta)\to(\Omega(M,\cH),\A)$ is a quasi-isomorphism, and
\[
\kappa=\sgn(\cK,\delta)\,(i_\infty)_*\circ\tau_{\cK,\delta},
\]
where $\kappa$ is composed with the isomorphism $\det\cK\cong\det H(M,\cH,\nabla)$ given by Hodge theory. Moreover $r(\delta)=\sum_{r\ge2}(\rk d_r|_{E_r^0}-\rk d_r|_{E_r^1})$.
\end{lemma}

\begin{proof}
The first two statements are the homological perturbation lemma for the special deformation retract $(i,\Pi_\cK,h)$ \cite[Lemma~2.4, \S3.2]{Cr}, with $h$ replaced by $-h$ to match the sign convention there. Filter $\cK$ by form degree. The spectral sequence of $(\cK,\delta)$ has $E_0=E_1=E_2=\cK$, and $i_\infty$ preserves the filtrations and induces the Hodge isomorphism $\cK\cong H(M,\cH,\nabla)$ on $E_2$. It therefore induces isomorphisms of all later pages and of $E_\infty$. By (KM3), $\sgn(\cK,\delta)\,\tau_{\cK,\delta}$ is the composite of the page isomorphisms, which correspond under $i_\infty$ to those defining $\kappa$. For the last statement, the sets $\delta b^0_r$ in the proof of (KM3) form a basis of $\delta(\cK^0)$, so $\rk\delta|_{\cK^0}=\sum_r\rk d_r|_{E_r^0}$, and likewise for odd parity. Since $i_\infty$ identifies the pages with $r\ge2$ and their differentials, these are the ranks in Theorem~\ref{thm:zero}.
\end{proof}

\subsection{Scaling}
As in \cite[\S2]{Mi}, for $t>0$ let $\rho_t=t^{N/2}$, where $N$ is the form-degree operator, and put
\[
D_t=t^{-1/2}\rho_t\A\rho_t^{-1}=\nabla+\sum_{j\ge2}t^{(j-1)/2}\A_{[j]}.
\]
This is a flat superconnection, and it extends continuously to $t=0$ with $D_0=\nabla$. Here we use $\partial=0$, because the family of \S1 is otherwise singular. Write
\[
w=\Str(N|_\cK)=\sum_k(-1)^kk\bigl(b_k(\cH^0)-b_k(\cH^1)\bigr),
\qquad \sdet(\rho_t|_\cK)=t^{w/2}.
\]

\begin{lemma}\label{lem:scaling}
For $t>0$, identify $\det H(M,\cH,D_t)$ with $\det H(M,\cH,\A)$ by $\rho_t^{-1}$.
\begin{enumerate}
\item $\|\cdot\|_{\RS(D_t)}=t^{\chi'(\A)/2}\|\cdot\|_{\RS(\A)}$.
\item Put $\delta_t=t^{-1/2}\rho_t\delta\rho_t^{-1}$ and $i^{(t)}_\infty=\rho_ti_\infty\rho_t^{-1}$. These are the operators of Lemma~\ref{lem:model} for $D_t$, and
\[
\Psi_t:=(\rho_t^{-1})_*\circ(i^{(t)}_\infty)_*\circ\tau_{\cK,\delta_t}=\sgn(\cK,\delta)\,t^{a}\kappa,\qquad a=-\tfrac12\bigl(w+r(\delta)\bigr).
\]
\end{enumerate}
\end{lemma}

\begin{proof}
(1) Although $\rho_t$ is not $\Omega(M)$-linear, it is an even bundle automorphism of $\Lambda T^*M\otimes\cH$. It is an isometry from $(\Omega(M,\cH),\A)$ with the metric $\langle\rho_t\cdot,\rho_t\cdot\rangle$ to $(\Omega(M,\cH),t^{1/2}D_t)$ with the fixed metric. Proposition~\ref{prop:metric} therefore identifies the Ray--Singer metrics of $\A$ and $t^{1/2}D_t$, as well as their derived Euler characteristics. Together with \eqref{eq:scale}, this gives $\chi'(D_t)=\chi'(\A)$ for $t>0$, extending \cite[Cor.~2.4]{Mi}, and proves (1).

(2) $\rho_t$ commutes with $i$ and $\Pi_\cK$, and $\rho_th\rho_t^{-1}=t^{-1/2}h$ because $h$ lowers form degree by one. Hence $h(D_t-\nabla)=\rho_t\,h\Atail\,\rho_t^{-1}$, which gives the formulas for $\delta_t$ and $i^{(t)}_\infty$. Since $\rho_t$ is an isomorphism $(\cK,\delta)\to(\cK,t^{1/2}\delta_t)$, (KM1) and (KM2) give
\[
\tau_{\cK,\delta_t}=t^{-r(\delta)/2}\,H(\rho_t)_*\,\tau_{\cK,\delta}\,(\det\rho_t|_\cK)^{-1}.
\]
The factor $t^{-r(\delta)/2}$ comes from scaling the differential, and $(\det\rho_t|_\cK)^{-1}$ acts by $t^{-w/2}$ on $\det\cK$. Finally $i^{(t)}_\infty\rho_t=\rho_ti_\infty$, so Lemma~\ref{lem:model} gives $\Psi_t=\sgn(\cK,\delta)\,t^{a}\kappa$.
\end{proof}

\subsection{Small eigenvalues}
For $t\ge0$ let $\Delta_t=D_tD_t^*+D_t^*D_t$, with adjoints taken in the fixed product metric. The difference $\Delta_t-\Delta_0$ has differential order at most one, and its coefficients are continuous in $t$ and vanish at $t=0$. By the elliptic estimate, $\|(\Delta_t-\Delta_0)u\|\le c(t)(\|\Delta_0u\|+\|u\|)$ with $c(t)\to0$, so $\Delta_t$ is self-adjoint on the Sobolev domain $H^2$ and $\Delta_t\to\Delta_0$ in the norm resolvent sense. Choose $\eigcut>0$ with $\operatorname{spec}\Delta_0\cap(0,2\eigcut)=\emptyset$. There is $t_0>0$ such that for $0\le t<t_0$ the number $\eigcut$ is not in the spectrum of $\Delta_t$ and the spectral projection $\Pi_t$ of $\Delta_t$ onto $[0,\eigcut)$ is continuous in $t$ \cite[IV-\S3]{Ka}, applied to the complexification. Its range $\Ksm$ has dimension $\dim\cK$, and $\cK_0=\cK$. Shrinking $t_0$, the polar decomposition
\[
U_t=\Pi_t\,(\Pi_\cK\Pi_t|_\cK)^{-1/2}\colon\cK\longrightarrow\Ksm
\]
is an orthogonal isomorphism that preserves parity, depends continuously on $t$ and satisfies $U_0=1$.

Since $D_t^2=0$, the operator $\Delta_t$ commutes with $D_t$ and $D_t^*$, so $\Ksm$ and $\Ksm^\perp$ are invariant under both. The restriction of $\Delta_t$ to $\Ksm^\perp$ has spectrum bounded below by $\eigcut$, so this subcomplex is acyclic. Its torsion contributes the scalar factor
\[
L_t=\exp\bigl(\tfrac12\zeta_{0,t}'(0)-\tfrac12\zeta_{1,t}'(0)\bigr),\qquad
\zeta_{i,t}(s)=\Tr\bigl(\Pi^{\co}_t(1-\Pi_t)\Delta_t^{-s}|_{\Omega(M,\cH)^i}\bigr),
\]
where $\Pi^{\co}_t$ is the orthogonal projection onto the closure of $\im D_t^*$. The remaining factor is the torsion of the finite complex $(\Ksm,D_t)$. Thus
\begin{equation}\label{eq:split}
\|\cdot\|_{\RS(D_t)}=L_t\,\|\cdot\|_{\Ksm,D_t},
\end{equation}
where $\|\cdot\|_{\Ksm,D_t}$ is the metric of Lemma~\ref{lem:finiteRS} on $\det H(\Ksm,D_t)=\det H(M,\cH,D_t)$. For $t=0$ the functions $\zeta_{i,0}$ are the zeta functions of Definition~\ref{def:RS} for $(\cH,\nabla)$, so $\|\cdot\|_{\RS(\nabla)}=L_0\,|\cdot|_{L^2}$.

\begin{lemma}\label{lem:cont}
$L_t\to L_0$ as $t\to0^+$.
\end{lemma}

\begin{proof}
Fix a parity $i$. By Lemma~\ref{lem:regular}(1), applied on $\Ksm^\perp$, we have
\[
\zeta_{i,t}(s)=\Gamma(s+1)^{-1}\int_0^\infty v^s\varphi_t(v)\,dv,
\]
initially for $\operatorname{Re}s$ large, where
\[
\begin{split}
\varphi_t(v)
&=\Tr\bigl(D_t^*D_te^{-v\Delta_t}(1-\Pi_t)|_{\Omega(M,\cH)^i}\bigr)\\
&=\Tr\bigl(D_t^*D_te^{-v\Delta_t}|_{\Omega(M,\cH)^i}\bigr)-\sum_\mu \mu e^{-v\mu}.
\end{split}
\]
Here $\mu$ runs over the eigenvalues of $D_t^*D_t$ on $\Ksm\cap\Omega(M,\cH)^i$. The finite-dimensional trace in the last term tends to zero as $t\to0$ and is uniformly bounded for $0<v\le1$ and $t$ near zero.

The coefficients of $\Delta_t$ and $D_t^*D_t$ are polynomial in $t^{1/2}$, and the principal symbol of $\Delta_t$ is that of the product metric. Thus $\Delta_t$ is a smooth family of generalized Laplacians in the parameter $t^{1/2}$, and heat kernel expansions with a parameter (see \cite[Lemma~2.49 and the proof of Thm.~2.48]{BGV}) give coefficients $a_j(t)$ continuous in $t$ and remainders uniform in $t$ for the first trace. Its coefficient $a_n(t)$ of $v^{-1}$ vanishes by the symbol-parity argument in Lemma~\ref{lem:regular}, since $n$ is odd. Consequently, after subtracting the terms with $j<n$, we can write
\[
\varphi_t(v)=\sum_{j=0}^{n-1}a_j(t)v^{(j-n-2)/2}+R_t(v),
\qquad |R_t(v)|\le c_0v^{-1/2}\quad(0<v\le1),
\]
with $c_0$ independent of $t$. For $v\ge1$, the uniform spectral gap on $\Ksm^\perp$ and a bound on the heat trace at $v=1$ uniform in $t$ give $|\varphi_t(v)|\le c_0e^{-\eigcut v}$, after increasing $c_0$. For fixed $v>0$, $\varphi_t(v)\to\varphi_0(v)$ as $t\to0$, because the kernel of $e^{-v\Delta_t}$ depends continuously on $t$ in the $C^\infty$ topology \cite[Thm.~2.48]{BGV} and the sum over $\mu$ tends to zero.

The meromorphic continuation near $s=0$ is therefore given by
\[
\begin{split}
\Gamma(s+1)\zeta_{i,t}(s)
={}&\sum_{j=0}^{n-1}\frac{a_j(t)}{s+(j-n)/2}
+\int_0^1 v^sR_t(v)\,dv\\
&+\int_1^\infty v^s\varphi_t(v)\,dv.
\end{split}
\]
None of the displayed denominators vanishes at $s=0$. On a fixed neighbourhood $|\operatorname{Re}s|<1/4$, the remainder bound and the exponential bound dominate both integrands and their $s$-derivatives uniformly in $t$. Dominated convergence, now applied to these subtracted integrals, proves continuity at $t=0$ of the value and the $s$-derivative at $s=0$ of $\Gamma(s+1)\zeta_{i,t}(s)$, and hence of $\zeta_{i,t}'(0)$. Applying this to both parities gives $L_t\to L_0$.
\end{proof}

\subsection{Proof of Theorem~\ref{thm:zero}}
Put $\tilde\delta_t=U_t^{-1}D_tU_t$ on $\cK$ and $\beta_t=U_t^{-1}\Pi_ti^{(t)}_\infty\colon(\cK,\delta_t)\to(\cK,\tilde\delta_t)$. This is a chain map, and $\beta_t\to1$ as $t\to0$ because $i^{(t)}_\infty\to i$, $\Pi_t\to\Pi_\cK$ and $U_t\to1$. Shrinking $t_0$, it is an isomorphism of complexes for $0<t<t_0$, with $\sdet\beta_t\to1$. Put
\[
\Phi_t=(\rho_t^{-1})_*\circ(j_t)_*\circ\tau_{\Ksm,D_t}\circ\det U_t\colon\det\cK\to\det H(M,\cH,\A),
\]
with $j_t\colon\Ksm\to\Omega(M,\cH)$ the inclusion. Since $\Ksm^\perp$ is acyclic, the projection $\Pi_t$ induces the inverse of $(j_t)_*$ on cohomology. Thus $(j_t)_*H(U_t)_*H(\beta_t)_*=(i_\infty^{(t)})_*$, and (KM1) gives $\Phi_t=\Psi_t\circ(\det\beta_t)^{-1}$.

Let $\ell\in\det\cK$ have $L^2$-norm one. Since $U_t$ is orthogonal, $\det U_t\,\ell$ has norm one in $\det\Ksm$, and Lemma~\ref{lem:finiteRS} gives $\|\tau_{\Ksm,D_t}(\det U_t\,\ell)\|_{\Ksm,D_t}=1$. The map $(\rho_t^{-1})_*\circ(j_t)_*$ sends this element of $\det H(\Ksm,D_t)$ to $\Phi_t\ell$. Moreover,
\[
\kappa\ell=\sgn(\cK,\delta)\,t^{-a}\Psi_t\ell
=\sgn(\cK,\delta)\,t^{-a}\sdet(\beta_t)\Phi_t\ell.
\]
Together with \eqref{eq:split}, Lemma~\ref{lem:scaling} and $\|\ell\|_{\RS(\nabla)}=L_0$, this gives
\[
\frac{\|\kappa\ell\|_{\RS(\A)}}{\|\ell\|_{\RS(\nabla)}}=t^{-\chi'(\A)/2-a}\,|\sdet\beta_t|\,\frac{L_t}{L_0}\qquad(0<t<t_0).
\]
The left side is a fixed positive number, while $|\sdet\beta_t|\,L_t/L_0\to1$ by Lemma~\ref{lem:cont}. A non-zero exponent of $t$ would force the right side to tend to zero or infinity. Hence $\chi'(\A)=-2a=w+r(\delta)$, and taking the limit shows that the left side equals one. With Lemma~\ref{lem:model}, this proves both the formula for $\chi'$ and the isometry. \qed

Theorem~B follows from Corollary~\ref{cor:reduction} and Theorem~\ref{thm:zero}. So does Theorem~C, since $\cH=H_\partial$ and the gauge transformations of \S\ref{sec:split} identify the pages $E_r$ with $r\ge1$, to which $\A_c$ contributes nothing.

\section{Reidemeister torsion: proof of Theorem~A}\label{sec:RT}

We recall the Reidemeister norm of Arias Abad and Sch\"atz \cite[\S3.1, Def.~3.8]{AS}. Write
\[
L_E=\det H(M,E,\A),\qquad L_{E^*}=\det H(M,E^*,\A^\vee),
\]
where $\A^\vee$ is the dual flat superconnection on the dual bundle $E^*$, not the adjoint \cite[Def.~2.6]{AS}. The canonical flat section of $\det(E\oplus E^*)\cong\RR$ gives a norm on the determinant line of the doubled simplicial complex \cite[Def.~3.6, Lemma~3.7]{AS}. Transport by the Knudsen--Mumford and de Rham isomorphisms gives a norm $\tau_{\mathrm{dbl}}$ on $L_E\otimes L_{E^*}$. Since $n$ is odd, the Poincar\'e pairing induces an isomorphism $\mathcal D_E\colon L_E\to L_{E^*}$, and the construction of \cite[Rem.~3.2, Def.~3.8]{AS} is
\[
\tau_R(v)=\tau_{\mathrm{dbl}}\bigl(v\otimes\mathcal D_Ev\bigr)^{1/2}.
\]
The square root, already present in \cite[Def.~3.8]{AS}, ensures that $\tau_R(zv)=|z|\tau_R(v)$. This construction follows the approach of Farber \cite[\S3.1]{AS} and, for an ordinary flat bundle, gives the Poincar\'e--Reidemeister metric of \cite[\S4.7]{Fa}.

The combinatorial comparison needed for Theorem~A is \cite[Cor.~3.12]{AS}. The following proposition restates it, makes explicit the comparison of the two spectral-sequence maps, and gives a direct verification of the Poincar\'e-duality compatibility in \cite[Prop.~3.11]{AS} using their duality lemmas.\footnote{Arias Abad and Sch\"atz consider isomorphisms of determinant lines only up to sign \cite[\S1]{AS}, and their natural isomorphism $\det C\cong\det H(C,\delta)$ \cite[Lemma~A.4]{AS} agrees with $\tau_{C,\delta}$ up to sign. We use the conventions of \S\ref{sec:conv} for both rows of \eqref{eq:square} and define $\mathrm{dR}$ by the actual induced map on $E_2$. Naturality then gives exact commutativity. No separate formula for the integration-map signs is needed.}

\begin{proposition}\label{prop:AS}
\begin{enumerate}
\item The integration map $\int$ of \cite[Cor.~2.14]{AS} preserves the filtrations, and on second pages it induces the classical de Rham isomorphism for $H_\partial$ up to signs determined by the grading conventions. With $\mathrm{dR}$ defined as the determinant of this induced map, the square \eqref{eq:square} commutes.
\item The isomorphism $\kappa_K$ in \eqref{eq:square} maps Farber's Poincar\'e--Reidemeister metric of $H_\partial$ isometrically to the Reidemeister norm $\tau_R$ of $(E,\A)$.
\end{enumerate}
\end{proposition}

\begin{proof}
(1) In a local trivialization write $\A=d+\sum_{j\ge0}a_{[j]}$, where $a_{[j]}$ is an endomorphism-valued $j$-form. The integration functor uses the modified $A_\infty$ de Rham map $\bar\psi$ of \cite[Def.~4.7]{AS2}. The map $\psi$ of \cite[Def.~3.23]{AS2} has source differential $-d$, and $\bar\psi=i\circ\psi\circ\tau$ changes the source and target conventions by degree-dependent signs. These changes preserve the vanishing property of \cite[Prop.~3.26(2)]{AS2}. The terms of the integration map are given by the components $\bar\psi_{m+1}$ applied to an $E$-valued form $\omega$ and $m$ of the forms $a_{[j]}$. For $m\ge1$, such a term vanishes if any input is a form of degree zero. Every non-zero term therefore has $j_i\ge1$ for the inserted forms $a_{[j_i]}$, and its cochain degree is
\[
\deg\omega+\sum_{i=1}^m j_i-m
=\deg\omega+\sum_{i=1}^m(j_i-1)\ge\deg\omega.
\]
The term $\bar\psi_1$ preserves form degree as well, so the integration map is filtered. Take local frames adapted to $\im\partial\subset\ker\partial$. By the gauge invariance of the integration map \cite[\S4.1]{AS2}, the terms preserving form degree on $E_1$ give the de Rham map with coefficients in $H_\partial$ and its induced flat connection. In local flat frames for $H_\partial$, this is ordinary integration with the degree-dependent signs from the source, target and suspension conventions of \cite[\S4.2]{AS2}. Thus its map on $E_2$ agrees with the classical de Rham isomorphism up to these signs, which act on determinant norms by factors of absolute value one. Define $\mathrm{dR}$ to be the determinant of this actual $E_2$-map. Since $\kappa$ and $\kappa_K$ are defined by the same conventions, the naturality of $\kappa$ (\S\ref{sec:conv}) shows that \eqref{eq:square} commutes.

(2) This is \cite[Cor.~3.12]{AS}. Arias Abad and Sch\"atz state it for the isomorphism of \cite[Rem.~3.10]{AS}, which comes from the form-degree filtration on forms. Their argument passes through the cellular cochains filtered by cochain degree in the proof of \cite[Thm.~3.9]{AS}, so it concerns $\kappa_K$, and part~(1) identifies the two maps. The proof also uses the compatibility of Poincar\'e duality with the spectral sequences \cite[Prop.~3.11]{AS}. We verify this compatibility directly from their Lemma~B.2, working with the induced pairings on the pages.

Write $H(E)=H(M,E,\A)$, and let $F_pH(E)$ consist of the classes with a representative in forms of degree at least $p$. Let $H(E^*)$ be the cohomology of the dual flat superconnection on the dual bundle $E^*$, with the corresponding filtration $\tilde F_pH(E^*)$. Let $\operatorname{gr}$ denote associated graded spaces, and let $\operatorname{Ann}$ denote the annihilator under the Poincar\'e pairing $H(E^*)\times H(E)\to\RR$, which is non-degenerate \cite[Rem.~2.7]{AS}. On the page $E_2=H(M,H_\partial)$ the induced pairing is Poincar\'e duality for the flat bundles $H_\partial$ and $(H_\partial)^*$, which is non-degenerate. By \cite[Lemma~B.2(4)]{AS}, applied from this page, the Poincar\'e pairing induces non-degenerate pairings
\[
\operatorname{gr}_pH(E^*)\times\operatorname{gr}_{n-p}H(E)\longrightarrow\RR
\]
on $E_\infty=\operatorname{gr}H$. The inclusion
\[
\tilde F_pH(E^*)\subseteq\operatorname{Ann}\bigl(F_{n-p+1}H(E)\bigr)
\]
is immediate, because representatives in these two filtration steps pair to a form of degree at least $n+1$. The non-degenerate pairings give
\[
\begin{split}
\dim\tilde F_pH(E^*)
&=\sum_{p'\ge p}\dim\operatorname{gr}_{p'}H(E^*)
=\sum_{q\le n-p}\dim\operatorname{gr}_{q}H(E)\\
&=\dim H(E)-\dim F_{n-p+1}H(E)
=\dim\operatorname{Ann}\bigl(F_{n-p+1}H(E)\bigr).
\end{split}
\]
The inclusion is therefore an equality, which is exactly the compatibility hypothesis of \cite[Lemma~B.2(5)]{AS}. No assertion about the form degrees of harmonic representatives is needed, and the rest of the proof of \cite[Cor.~3.12]{AS} applies unchanged.
\end{proof}

\begin{proof}[Proof of Theorem~A]
By Proposition~\ref{prop:AS}, $\kappa$ is an isometry for the Reidemeister norms of $H_\partial$ and of $(E,\A)$, transported to de Rham cohomology by the de Rham isomorphisms. For the flat bundle $H_\partial$, the construction of $\tau_R$ \cite[\S3.1, Def.~3.8]{AS} is Farber's Poincar\'e--Reidemeister metric \cite[\S4.7]{Fa}, which is defined over $\RR$ or $\CC$ \cite[\S1.5]{Fa}. Multiplicativity of the Reidemeister norm is \cite[Lemma~3.16(2)]{AS}. Applying it to $H_\partial=H_\partial^0\oplus H_\partial^1$, with the norm on the odd summand taken on the dual determinant line, reduces the comparison to ungraded bundles. The Ray--Singer metric has the same direct-sum and parity-reversal properties. On odd-dimensional manifolds the two metrics are equal for arbitrary flat bundles and arbitrary metrics \cite[Thm.~6.2]{Fa}. This is a consequence of \cite[Thm.~0.2]{BZ}, which is stated for complex flat bundles and applies to the real flat bundle $H_\partial$ through its complexification, by Lemma~\ref{lem:cplx}. Theorem~B completes the proof.
\end{proof}

\begin{remark}[Complex coefficients]\label{rem:complex}
Let $E$ be a $\ZZ_2$-graded complex vector bundle with a Hermitian metric, and let $\A$ be a complex-linear flat superconnection. The definitions of \S\ref{sec:conv} and the proofs in \S\ref{sec:metric} and \S\S\ref{sec:split}--\ref{sec:zero} apply verbatim over $\CC$, so Theorems~B and~C hold for $(E,\A)$. Goette's construction and Ma's theorem are available over $\CC$, so the first statement of Corollary~\ref{cor:morse} holds for a complex flat bundle $F$ with a Hermitian metric. The construction of $\tau_R$ in \cite{AS} is over $\RR$. Its complex version uses complex determinant lines, the complex-bilinear Poincar\'e pairing and the canonical section of $\det_\CC(E\oplus E^*)$. Farber's metric is defined over $\CC$ as well, and the arguments of \cite{AS} and of Proposition~\ref{prop:AS} hold over either ground field. With this construction Theorem~A over $\CC$ follows from Theorem~B over $\CC$ and \cite[Thm.~6.2]{Fa}, and so does the second statement of Corollary~\ref{cor:morse}.
\end{remark}

\section{Morse flat superconnections}\label{sec:morse}

A main source of geometric examples of flat superconnections with non-zero higher components is the work of Goette on Morse theory and higher torsion invariants of families \cite{Go1,Go2}. The aim of this work is to compare the analytic torsion forms of Bismut and Lott \cite{BL} with finite-dimensional torsion forms built from family Morse data, and, under the additional hypotheses needed for their definition and comparison, with the Igusa--Klein higher Franz--Reidemeister torsion. These higher invariants depend on the smooth structure of the bundle. In particular, Goette uses higher analytic torsion to distinguish bundles with diffeomorphic fibres that are homeomorphic but not diffeomorphic as bundles \cite[\S0]{Go1}.

\subsection{The Morse model}\label{ssec:goette-model}
Let $p\colon X\to B$ be a smooth fibre bundle with closed fibres $Z_b=p^{-1}(b)$ over a closed oriented manifold $B$, and let $F\to X$ be a flat real vector bundle with connection $\nabla^F$ and metric $h^F$. Let $\Mf\colon X\to\RR$ be a function whose restriction to every fibre is Morse, and choose a horizontal distribution, a vertical metric and an auxiliary metric on $F$ satisfying the normal form assumptions of \cite[Prop.~1.10]{Go1}. Goette's construction is available over $\RR$ or $\CC$ \cite[\S5.2]{Go3}. We use three of its properties \cite[Def.~1.1, Def.~1.60, Thm.~1.61]{Go1}.
\begin{enumerate}
\item The critical points define a finite-rank graded bundle $V=\bigoplus_qV^q$ with
\[
V_b^q=\bigoplus_{y}F_y\otimes o(T_y^uZ_b),
\]
the sum being over the critical points $y$ of $\Mf|_{Z_b}$ of index $q$. Parallel transport in $F$ along the critical set and the orientation local systems give a flat connection $\nabla^V$.
\item There is a flat superconnection
\begin{equation}\label{eq:goette-superconnection}
\A'=\nabla^V+\sum_{j\ge0}a'_j,\qquad a'_j\in\Omega^j(B,\End^{1-j}V),
\end{equation}
where the superscript on $\End V$ is the change in Morse index. With the grading of $V$ by the parity of $q$, this is a flat superconnection in the sense of \S1. The operator $\A'-\nabla^V$ strictly increases critical values: it maps the summand at $y$ only to summands at critical points with larger value of $\Mf$. The zero-form part $\partial=a'_0$ computes the fibre cohomology. At each $b$ it is conjugate, by $1+c$ with $c\in\End^0V$ strictly increasing critical values, to the Thom--Smale differential of a Smale gradient on a nearby fibre \cite[(1.53)]{Go1}.
\item There is a cochain map
\begin{equation}\label{eq:morse-integration}
I\colon(\Omega(X,F),d^F)\longrightarrow(\Omega(B,V),\A')
\end{equation}
that is linear over $\Omega(B)$, maps the filtration of $\Omega(X,F)$ by horizontal degree \cite[Def.~1.4]{Go1} to the filtration of $\Omega(B,V)$ by form degree, and induces isomorphisms of spectral sequences from $E_1$ onwards. Under the local triangular identifications with ordinary Thom--Smale complexes, its map on $E_1$ agrees with that induced by integration over unstable cells. More precisely, the proof of \cite[Thm.~1.61, (1.62)]{Go1} compares $I$ with the local unstable-cell integration maps by a triangular conjugation and a cochain homotopy.
\end{enumerate}
A gradient that is Morse--Smale on every fibre is not required. Across parameter values where the Smale condition fails, the local Morse complexes are compared by continuation maps, whose compatibility is expressed by higher homotopies. The forms $a'_j$ retain this family data in a single flat superconnection \cite[\S1]{Go1}. The spectral sequence of the filtration by horizontal degree is the Leray spectral sequence \cite[\S2]{Ma}, so $I$ identifies $(H(V,\partial),\nabla^H)$ with the flat bundle $H(Z,F)$ and $H(B,V,\A')$ with $H(X,F)$.

\begin{corollary}\label{cor:morse}
Suppose that $B$ is odd-dimensional and that the fibre $Z$ is even-dimensional. Then
\[
I_*\colon \det H(X,F)\longrightarrow \det H(B,V,\A')
\]
is an isometry from the Ray--Singer metric of $(X,F)$ to that of $(V,\A')$. Under the de Rham isomorphism of Theorem~A for $(V,\A')$, both equal the Reidemeister norm $\tau_R$ of $(V,\A')$ on a triangulation of $B$.
\end{corollary}

This is a statement about a norm on the single determinant line $\det H(B,V,\A')$. It does not concern Goette's torsion form on $B$, which is compared with the analytic torsion form of Bismut and Lott \cite{BL} in \cite[Thm.~0.1]{Go1}.

\begin{proof}
By property (3), $I_*$ intertwines the isomorphism $\kappa_{V}$ of Theorem~B for $(V,\A')$ with the isomorphism $\kappa_{\mathrm{Leray}}\colon\det H(B,H(Z,F))\to\det H(X,F)$ given by the Leray spectral sequence. Ma's functoriality theorem holds for arbitrary complex flat bundles and Hermitian metrics. In the form \cite[(0.6)]{Ma}, applied to $F\otimes\CC$ and combined with Lemma~\ref{lem:cplx}, it says that the logarithm of the ratio of the Ray--Singer metrics under $\kappa_{\mathrm{Leray}}$ is
\[
\int_Be(TB)\,\log T(Z,h^F)-\tfrac12\int_X\tilde e(TX)\,\theta(F,h^F).
\]
Here $e(TB)$ is the Euler form of $TB$, $T(Z,h^F)$ is the Ray--Singer torsion of the fibres, a function on $B$, $\tilde e(TX)$ is the Chern--Simons form that transgresses between the Euler forms of $TX$ for two connections, and $\theta(F,h^F)=\Tr\bigl((h^F)^{-1}\nabla^Fh^F\bigr)$ is the Kamber--Tondeur form of $F$. The first term vanishes because $\dim B$ is odd. The second term vanishes because $\dim X$ is odd, since the Euler form of a bundle of odd rank and its Chern--Simons transgression vanish identically. Hence $\kappa_{\mathrm{Leray}}$ is an isometry. For unimodular $F$ this is due to L\"uck, Schick and Thielmann \cite{LST}. Theorem~B says that $\kappa_V$ is an isometry. The Ray--Singer metric of the flat bundle $H(Z,F)$ on the odd-dimensional $B$ does not depend on its metric (Proposition~\ref{prop:metric}, or \cite[Thm.~0.1]{BZ}), so $I$ induces an isometry between the Ray--Singer metrics of $H(Z,F)$ and $H(V,\partial)$, and $I_*$ is an isometry. The last statement is Theorem~A for $(V,\A')$.
\end{proof}

\begin{remark}
If $\dim Z$ is odd, the same argument applies whenever the Chern--Simons term of \cite[(0.6)]{Ma} vanishes, for instance when $h^F$ has a parallel volume form.
\end{remark}

\subsection{A double cover of a mapping torus}\label{ssec:doublecover}

We give an example with both a non-zero Morse differential and a non-zero higher term.

\emph{The fibration and the Morse function.}
Identify $S^3$ with the unit quaternions and put $\Psi(b,v)=(b,bv)$ on $S^3\times S^3$. For $m\ge1$ let
\[
X_m=\bigl([0,2\pi m]\times S^3\times S^3\bigr)\big/\bigl((2\pi m,b,v)\sim(0,b,b^mv)\bigr).
\]
Projection onto $b$ gives a bundle $X_m\to S^3$. Its fibre is the mapping torus of left multiplication by $b^m$, which is isotopic to the identity, so the fibre is $S^1\times S^3$. We take trivial coefficients. There is a connected $m$-fold covering $r_m\colon X_m\to X_1$ over $S^3$, given on the $j$-th interval by
\[
r_m[\phi,b,v]=[\phi-2\pi j,b,b^jv],
\qquad 2\pi j\le\phi\le2\pi(j+1),\quad 0\le j<m.
\]
To construct a fibrewise Morse function on $X_1$, fix $0<\gamma<\pi/2$ and a smooth cutoff $\mu\colon\RR\to[0,1]$ with $\mu=0$ on $(-\infty,\pi+\gamma]$ and $\mu=1$ on $[2\pi-\gamma,\infty)$. Set
\[
\xi(b,\phi)=1-\mu(\phi)+\mu(\phi)b^{-1},\qquad
\Mf_1(\phi,b,v)=\cos\phi+\varepsilon\langle v,\xi(b,\phi)\rangle,
\]
where $\langle\cdot,\cdot\rangle$ is the Euclidean inner product on the quaternions. The identity $\langle bv,1\rangle=\langle v,b^{-1}\rangle$ shows that $\Mf_1$ descends to $X_1$. Choose $0<\varepsilon<1$ with $\varepsilon\sup|\partial_\phi\xi|<\sin\gamma$. Then $\partial_\phi\Mf_1\ne0$ on $[\pi+\gamma,2\pi-\gamma]$, and the critical points on every fibre are
\[
\begin{array}{c|c|c}
(\phi,v)&\text{Morse index}&\Mf_1\\\hline
(\pi,-1)&0&-1-\varepsilon\\
(0,-1)&1&1-\varepsilon\\
(\pi,+1)&3&-1+\varepsilon\\
(0,+1)&4&1+\varepsilon
\end{array}
\]
Put $\Mf_m=r_m^*\Mf_1$, with the vertical gradient in the normal form of \cite[Prop.~1.10]{Go1}. We focus on $m=2$. Every fibre has two critical points in each index $q\in\{0,1,3,4\}$, so $V^q=\langle\Mgen{q}{0},\Mgen{q}{1}\rangle$.

\emph{The Morse differential.}
Choose orientations so that, in each sphere degree, the degree-zero differential is the cellular coboundary for a circle with two vertices and two edges:
\begin{equation}\label{eq:doublecover-partial}
\partial|_{V^0\to V^1}=\partial|_{V^3\to V^4}
=\begin{pmatrix}-1&1\\1&-1\end{pmatrix},
\qquad \partial|_{V^1\oplus V^4}=0.
\end{equation}
In particular, the integral fibre cohomology classes
\[
u_3=[\Mgen{3}{0}+\Mgen{3}{1}]\in H^3(V,\partial),\qquad
u_1=[\Mgen{1}{0}]=[\Mgen{1}{1}]\in H^1(V,\partial)
\]
are primitive generators.

\emph{Integral cohomology.}
Let $\alpha,\beta\in H^3(S^3\times S^3;\ZZ)$ be the generators from the two factors. Quaternion multiplication pulls the generator back to the sum of the generators, so
\[
(\Psi^2)^*\alpha=\alpha,\qquad (\Psi^2)^*\beta=\beta+2\alpha.
\]
The action is the identity on $H^0$ and $H^6$. Write $R_q=((\Psi^2)^*-1)|_{H^q(S^3\times S^3;\ZZ)}$. The Wang sequence gives
\[
0\longrightarrow\operatorname{coker}R_{q-1}
\longrightarrow H^q(X_2;\ZZ)
\longrightarrow\ker R_q\longrightarrow0.
\]
On $H^3$, the map sends $\beta$ to $2\alpha$ and $\alpha$ to zero. Hence
\[
H^q(X_2;\ZZ)=
\begin{cases}
\ZZ,&q=0,1,3,6,7,\\
\ZZ\oplus\ZZ/2,&q=4,\\
0,&\text{otherwise}.
\end{cases}
\]
Since the base $S^3$ is simply connected, the integral Leray page $E_2^{p,q}=H^p(S^3;H^q(S^1\times S^3;\ZZ))$ is $\ZZ$ for $p\in\{0,3\}$ and $q\in\{0,1,3,4\}$, and zero otherwise. Its only possible non-zero differential is $d_3\colon E_3^{0,3}\to E_3^{3,1}$. The groups just computed force this map to be injective with cokernel $\ZZ/2$. After choosing orientations, it is multiplication by $2$.

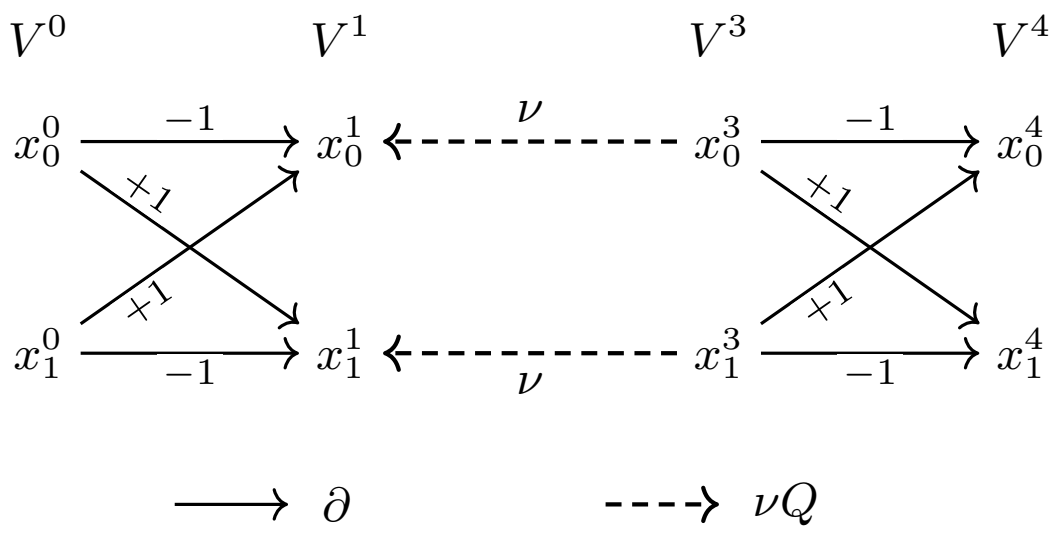
\begin{figure}[tbp]
\centering
\begin{tikzpicture}[
  x=1cm,y=1cm,
  every node/.style={font=\small},
  generator/.style={inner sep=3pt},
  differential/.style={->,semithick},
  higher/.style={->,thick,densely dashed},
  coefficient/.style={font=\scriptsize,fill=white,inner sep=1pt}
]
\node at (0,3.5) {$V^0$};
\node at (2,3.5) {$V^1$};
\node at (4.5,3.5) {$V^3$};
\node at (6.5,3.5) {$V^4$};
\node[generator] (a0) at (0,2.8) {$\Mgen{0}{0}$};
\node[generator] (a1) at (0,1.4) {$\Mgen{0}{1}$};
\node[generator] (b0) at (2,2.8) {$\Mgen{1}{0}$};
\node[generator] (b1) at (2,1.4) {$\Mgen{1}{1}$};
\node[generator] (c0) at (4.5,2.8) {$\Mgen{3}{0}$};
\node[generator] (c1) at (4.5,1.4) {$\Mgen{3}{1}$};
\node[generator] (d0) at (6.5,2.8) {$\Mgen{4}{0}$};
\node[generator] (d1) at (6.5,1.4) {$\Mgen{4}{1}$};
\draw[differential] (a0) -- node[coefficient,above] {$-1$} (b0);
\draw[differential] (a1) -- node[coefficient,below] {$-1$} (b1);
\draw[differential] (a0) -- node[coefficient,pos=.25,above,sloped] {$+1$} (b1);
\draw[differential] (a1) -- node[coefficient,pos=.25,below,sloped] {$+1$} (b0);
\draw[differential] (c0) -- node[coefficient,above] {$-1$} (d0);
\draw[differential] (c1) -- node[coefficient,below] {$-1$} (d1);
\draw[differential] (c0) -- node[coefficient,pos=.25,above,sloped] {$+1$} (d1);
\draw[differential] (c1) -- node[coefficient,pos=.25,below,sloped] {$+1$} (d0);
\draw[higher] (c0) -- node[above] {$\nu$} (b0);
\draw[higher] (c1) -- node[below] {$\nu$} (b1);
\draw[differential] (.9,.4) -- (1.65,.4);
\node[anchor=west] at (1.75,.4) {$\partial$};
\draw[higher] (3.75,.4) -- (4.5,.4);
\node[anchor=west] at (4.6,.4) {$\nu Q$};
\end{tikzpicture}
\caption{The model \eqref{eq:doublecover-model}, with the operators $\partial$ (solid) and $\nu Q$ (dashed).}
\label{fig:doublecover}
\end{figure}

\emph{The higher term.}
Let $\nu$ be a $3$-form on $S^3$ with $\int_{S^3}\nu=1$. Up to a gauge transformation with degree-zero part the identity, the Morse superconnection is
\begin{equation}\label{eq:doublecover-model}
\A_2=d+\partial+\nu Q,\qquad Q(\Mgen{3}{i})=\Mgen{1}{i}\quad(i=0,1),
\qquad Q|_{V^0\oplus V^1\oplus V^4}=0.
\end{equation}
Figure~\ref{fig:doublecover} shows this model. To see this form, use property (2). The critical sections have trivial orientation local systems over $S^3$, so $\nabla^V=d$. The critical values are equal within each index, so the term of form degree one vanishes. A term of form degree two would lower the critical value and is therefore also excluded. Thus the superconnection has the form $d+\partial+\Theta$, with $\Theta\in\Omega^3(S^3,\Hom(V^3,V^1))$ and $\partial$ constant by flatness. On the product fibre at $b=1$, its value is \eqref{eq:doublecover-partial}. Critical values coincide within each index, so the conjugation in property (2) is trivial and $\partial$ is the Thom--Smale differential of any Smale gradient on a single fibre.

Put $\bar\Theta=\int_{S^3}\Theta$. A gauge transformation of form degree two replaces $\Theta$ by $\nu\bar\Theta$. The topological input is $d_3=\pm2$, which after choosing orientations says that $\bar\Theta$ induces $u_3\mapsto2u_1$ on fibre cohomology. The map $Q$ induces the same map, since
\begin{equation}\label{eq:doublecover-transgression}
Q_*u_3=[\Mgen{1}{0}+\Mgen{1}{1}]=2u_1,\qquad d_3(u_3)=2[\nu]\otimes u_1.
\end{equation}
Hence $\bar\Theta-Q$ induces zero on fibre cohomology and is $[\partial,\cdot\,]$-exact. A gauge transformation of form degree three removes $\nu(\bar\Theta-Q)$ and gives \eqref{eq:doublecover-model}.

\emph{The determinant norm.}
Let $\ell_\ZZ(E_2)$ be the determinant generator defined by integral bases on the Leray page, and let $\ell_\ZZ(H)$ be the generator from the surviving free classes. Both are defined up to sign. The cancelled class $u_3$ has odd total degree $3$, while $[\nu]\otimes u_1$ has even total degree $4$, so
\[
\kappa\bigl(\ell_\ZZ(E_2)\bigr)=\pm\tfrac12\ell_\ZZ(H).
\]
The Gauss--Manin bundle consists of four trivial graded lines. Classical Cheeger--M\"uller for $S^3$, whose integral cohomology is torsion free, gives $\|\ell_\ZZ(E_2)\|_{\RS}=1$. Theorem~B therefore yields
\[
\|\ell_\ZZ(H)\|_{\RS(S^3,V,\A_2)}=2.
\]
Independently, the Reidemeister norm of the determinant generator of integral cohomology modulo torsion is $\prod_q|\operatorname{tors}H^q|^{(-1)^q}$. Since $X_2$ has just one torsion group, namely $\ZZ/2$ in degree four, classical Cheeger--M\"uller gives $\|\ell_\ZZ(H(X_2))\|_{\RS(X_2)}=2$. Here the triangular degree-zero conjugations in property (3) are the identity, because critical values coincide within each index. Thus $I$ induces the integral unstable-cell identifications on fibre cohomology. The subsequent gauge transformations have degree-zero part the identity, so they preserve this identification on $E_1$. Each non-zero free cohomology group of $X_2$ has just one surviving free Leray term, so the determinant generators agree under $I_*$, up to sign. This is consistent with Corollary~\ref{cor:morse}. Theorem~C also gives $\chi'(\A_2)=-1$, since the graded Betti contribution is zero and $d_3$ has rank one from odd to even parity.

The model also shows what the example tests. For the metric on $V$ in which the generators are orthonormal, $\cH$ is spanned by the sums $\Mgen{q}{0}+\Mgen{q}{1}$, the bundle $E_c$ is spanned by the differences, and $Q$ preserves both. Hence $\A_2=\A_\cH\oplus\A_c$ already. On $\cH$ the higher term sends the unit vector $\Mhat{3}=(\Mgen{3}{0}+\Mgen{3}{1})/\sqrt2$ to $\nu\Mhat{1}$, with $\Mhat{1}=(\Mgen{1}{0}+\Mgen{1}{1})/\sqrt2$. This has the form of a unit top-degree flux, but it maps one line of $\cH$ to another, so it is not the scalar flux of \cite[Prop.~5.1]{MW1} and \cite[Thm.~3.7]{Mi}. The summand $\A_c$ contains $\nu Q$, which a gauge transformation of form degree three removes, and it contributes nothing by Propositions~\ref{prop:normal} and~\ref{prop:double}. The factor $2$ is a lattice index: $u_3=\sqrt2\,\Mhat{3}$ and $u_1=\Mhat{1}/\sqrt2$, so $Q_*u_3=2u_1$. The example therefore illustrates Theorem~B with a non-zero Morse differential, in a case where the flux computation is elementary. In the same way, the $m$-fold cover $X_m$ gives $d_3=m$ and norm $m$.

\newcommand{\Nk}{\mathrm{Nil}_k}
\newcommand{\Xkm}{X_{k,m}}
\newcommand{\Akm}{\A_{k,m}}
\newcommand{\Tx}{\mathsf{x}}
\newcommand{\Ty}{\mathsf{y}}
\newcommand{\Tu}{\mathsf{u}}
\newcommand{\Tgen}[1]{\mathsf{e}_{#1}}
\newcommand{\Tgenp}{\mathsf{e}'_1}
\newcommand{\Ttw}{\mathsf{T}}
\newcommand{\Lk}{\mathcal{L}_k}
\newcommand{\Wk}{\mathcal{W}}
\newcommand{\nuB}{\nu_B}
\newcommand{\areaS}{\omega_{S^2}}
\newcommand{\gkm}{\gcd(k,m)}
\newcommand{\lkm}{\operatorname{lcm}(k,m)}

\subsection{Gauss--Manin monodromy}\label{ssec:monodromy}

In \S\ref{ssec:doublecover} the Gauss--Manin connection is trivial and $a'_1=0$. We now add the monodromy of a Dehn twist. This produces a summand on which $d_3$ vanishes although the higher term cannot be removed by a gauge transformation. Fix integers $k,m\ge1$.

\emph{Morse functions on the torus.}
Let $T^2=\RR^2/\ZZ^2$ have coordinates $(\Tx,\Ty)$ and let
\[
\Nk=(\RR\times T^2)\big/\bigl((\vartheta+1,\Tx,\Ty)\sim(\vartheta,\Tx+k\Ty,\Ty)\bigr)
\]
be the mapping torus of the $k$-th power of a Dehn twist, a bundle over $S^1=\RR/\ZZ$ with coordinate $\vartheta$. Take $\Ty\in[-\tfrac14,\tfrac34)$ and put $\Tu=\Tx+k\vartheta\Ty$. Let $\varrho$ be a smooth cutoff on $[-\tfrac14,\tfrac34]$ equal to $1$ near $0$ and $\tfrac12$ and equal to $0$ near the endpoints. For small $\varepsilon>0$, the function and the vertical metric
\[
h_k=\cos2\pi\Ty+\varepsilon\,\varrho(\Ty)\cos2\pi\Tu,\qquad
d\Tu^2+d\Ty^2=(d\Tx+k\vartheta\,d\Ty)^2+d\Ty^2
\]
descend to $\Nk$. In the coordinates $(\Tu,\Ty)$ both are independent of $\vartheta$, and $\vartheta$ enters only through the gluing across $\Ty=-\tfrac14$, where $\Tu$ jumps by $k\vartheta$. The critical points on every fibre are
\[
\begin{array}{c|c|c|c}
\text{generator}&(\Tu,\Ty)&\text{Morse index}&h_k\\\hline
\Tgen0&(\tfrac12,\tfrac12)&0&-1-\varepsilon\\
\Tgen1&(0,\tfrac12)&1&-1+\varepsilon\\
\Tgenp&(\tfrac12,0)&1&1-\varepsilon\\
\Tgen2&(0,0)&2&1+\varepsilon
\end{array}
\]
The lines $\Tu=\tfrac12$ are invariant under the descending gradient flow. From $\Tgenp$, the arc with $\Ty$ increasing ends at $\Tgen0$. The arc with $\Ty$ decreasing crosses $\Ty=-\tfrac14$ and continues with $\Tu=\tfrac12+k\vartheta$, so it ends at $\Tgen1$ if $k\vartheta\equiv\tfrac12$ modulo $1$ and at $\Tgen0$ otherwise. Hence the Smale condition fails exactly on the $k$ fibres $\vartheta=(j+\tfrac12)/k$ with $0\le j<k$. Each of them has one flow line from $\Tgenp$ to $\Tgen1$, a handle slide, and all of these have the same sign because $\tfrac12+k\vartheta$ has positive derivative.

The critical sections have trivial orientation local systems, and the Morse differential of $T^2$ vanishes. The only map between generators that preserves the index and increases critical values is $\Tgen1\mapsto\Tgenp$. So by property (2) the Morse superconnection of $\Nk\to S^1$ is $d+\omega\Ttw$, where $\omega$ is a $1$-form on $S^1$, $\Ttw(\Tgen1)=\Tgenp$, and $\Ttw$ vanishes on the other generators. Normalize the frame at $\vartheta=0$ using the integral unstable-cell classes $[d\Tx]$ and $[d\Ty]$. The pullback of the gluing diffeomorphism sends $[d\Tx]$ to $[d\Tx]+k[d\Ty]$ and fixes $[d\Ty]$. Thus $\int_{S^1}\omega=k$ after choosing orientations, and a gauge transformation equal to the identity at $\vartheta=0$ replaces $\omega$ by $k\,d\vartheta$. This specifies the integral local system inside the resulting flat bundle.

\emph{The fibration and the model.}
Let $B=S^1\times S^2$, let $\mathrm{pr}_1\colon B\to S^1$ be the projection, and let $\varpi\colon B\to S^3$ have degree one. With $X_m$ as in \S\ref{ssec:doublecover}, put
\[
\Xkm=\mathrm{pr}_1^*\Nk\times_B\varpi^*X_m.
\]
This is a bundle over $B$ with fibre $Z=T^2\times S^1\times S^3$, and we take trivial coefficients. The fibrewise Morse function is the sum of the pullbacks of $h_k$ and $\Mf_m$, so the critical points and their Morse indices are products.

The argument of \S\ref{ssec:doublecover} applies to every $m$. On $V_m=\langle\Mgen{q}{i}:q\in\{0,1,3,4\},\ i\in\ZZ/m\rangle$ it gives the model $d+\partial+\nu Q$ of $X_m$, with
\[
\partial\Mgen{0}{i}=\Mgen{1}{i-1}-\Mgen{1}{i},\qquad
\partial\Mgen{3}{i}=\Mgen{4}{i-1}-\Mgen{4}{i},\qquad
Q\Mgen{3}{i}=\Mgen{1}{i},
\]
and with $\partial$ and $Q$ zero on the other generators. The fibre cohomology has the integral generators $u_0=[\sum_i\Mgen{0}{i}]$, $u_1=[\Mgen{1}{0}]$, $u_3=[\sum_i\Mgen{3}{i}]$ and $u_4=[\Mgen{4}{0}]$, and $Q_*u_3=mu_1$. Let $\areaS$ be the area form of $S^2$ with total area one, and let $\nuB=d\vartheta\wedge\areaS$. Since $\varpi$ has degree one, a gauge transformation of form degree two turns the tensor product of the two superconnections into
\begin{equation}\label{eq:monodromy-model}
\Akm=d+\partial+k\,d\vartheta\,\Ttw+\nuB Q
\qquad\text{on}\qquad V=\langle\Tgen0,\Tgen1,\Tgenp,\Tgen2\rangle\otimes V_m.
\end{equation}
Here $\Ttw$ acts on the first factor, and $\partial$ and $Q$ act on the second, with the Koszul sign for $\partial$. This model satisfies properties (1)--(3), with $I$ the composite of the integration map of $\Xkm\to\varpi^*X_m$, which is pulled back from $\Nk\to S^1$, and the integration map of $\varpi^*X_m\to B$, which extends to these coefficients because it is $\Omega(B)$-linear. The proof of Corollary~\ref{cor:morse} uses only these properties, so the corollary applies. In \eqref{eq:monodromy-model}, $\nabla^V=d$ and $a'_1=k\,d\vartheta\,\Ttw\ne0$.

\emph{The spectral sequence.}
Let $\Lk=\langle\Tgen1,\Tgenp\rangle$ with the flat connection $\nabla=d+k\,d\vartheta\,\Ttw$, whose monodromy is unipotent and not orthogonal for any metric. Then
\[
H_\partial=\bigl(\langle\Tgen0\rangle\oplus\Lk\oplus\langle\Tgen2\rangle\bigr)\otimes\langle u_0,u_1,u_3,u_4\rangle.
\]
Since $\nabla\Tgen1=k\,d\vartheta\,\Tgenp$ and $\nabla(\areaS\Tgen1)=k\,\nuB\Tgenp$, the real cohomology of $\Lk$ over $B$ is generated by $\Tgenp$, $d\vartheta\,\Tgen1$, $\areaS\Tgenp$ and $\nuB\Tgen1$ in degrees $0,1,2,3$, and over $\ZZ$ the degrees $1$ and $3$ have additional torsion $\ZZ/k$, generated by $d\vartheta\,\Tgenp$ and $\nuB\Tgenp$. By \eqref{eq:monodromy-model} the only non-zero differential on $E_2$ is $d_3$, which sends $c\,u_3$ to $m\nuB c\,u_1$ for $c\in\{\Tgen0,\Tgenp,\Tgen2\}$. For $c=\Tgen0,\Tgen2$ this is multiplication by $m$, from odd to even total degree. For $c=\Tgenp$ it vanishes over $\RR$.

\begin{proposition}\label{prop:monodromy}
Let $\Mhat{3}=m^{-1/2}\sum_i\Mgen{3}{i}$ and $\Mhat{1}=m^{-1/2}\sum_i\Mgen{1}{i}$. Then $\Wk=\Lk\otimes\langle\Mhat{3},\Mhat{1}\rangle$ is a direct summand of $(V,\Akm)$, on which the superconnection is $\nabla+\nuB Q$. On $H(B,\Wk,\nabla)$ the differential $d_3$ vanishes, and $\nabla+\nuB Q$ is not gauge equivalent to $\nabla$.
\end{proposition}

\begin{proof}
The operators $\partial$ and $Q$ preserve the splitting of $V_m$ into sums and differences of generators with equal $q$, and $\Ttw$ preserves $\Lk$ and $\langle\Tgen0,\Tgen2\rangle$. This gives the summand, and $d_3$ vanishes because $\nuB\Tgenp=k^{-1}\nabla(\areaS\Tgen1)$. All of $\Wk$ has even parity, so a gauge transformation is $g=g_0+g_2$ with $g_j\in\Omega^j(B,\End\Wk)$ and $g_0$ invertible. Suppose that $g(\nabla+\nuB Q)=\nabla g$. In form degree one, $g_0$ is parallel, so $g_0=1\otimes G_1+\Ttw\otimes G_2$ with $G_1$ invertible. In form degree three, $g_0\nuB Q=[\nabla,g_2]$. The partial trace over $\Lk$ annihilates commutators with $\Ttw$, and it sends $g_0Q$ to $2G_1Q$ since $\Tr\Ttw=0$. So $2\nuB G_1Q$ is the exterior derivative of the partial trace of $g_2$. Integrating over $B$ gives $G_1Q=0$, which contradicts $Q\ne0$.
\end{proof}

The higher term on $\Wk$ is therefore invisible to the pages and, by Theorem~B, to the Ray--Singer metric, although it changes the integral cohomology computed below. The superconnections $\nabla+c\,\nuB Q$ with $c>0$ are gauge equivalent by constant automorphisms of $\langle\Mhat{3},\Mhat{1}\rangle$, and their limit $\nabla$ as $c\to0$ is not, as in the rescaling of \S\ref{sec:zero}.

\emph{Integral cohomology.}
We construct an integral mapping-torus model to fix both the torsion groups and the cohomology lattices. Regard $\varpi^*X_m$ as a sphere bundle over
\[
C=S^1_\vartheta\times S^2\times S^1_\phi.
\]
It is the unit sphere bundle of a quaternionic line bundle, obtained by pulling back the clutching map $b\mapsto b^m$ of $X_m$. Put $\omega=[\areaS]$ and let $\eta$ generate $H^1(S^1_\phi;\ZZ)$. Its Euler class is $m[d\vartheta]\omega\eta$. Quaternionic line bundles over a four-dimensional CW complex are classified by this class, since $BS^3$ is $3$-connected and $\pi_4(BS^3)=\ZZ$. We can therefore cut the $\vartheta$-circle and use a clutching map $a\colon S^2\times S^1_\phi\to S^3$ of degree $m$. Choose $a$ constant on $S^2\vee S^1_\phi$. It follows that $\Xkm$ is the mapping torus of the diffeomorphism
\[
\Phi\colon Y=T^2\times S^2\times S^1_\phi\times S^3\longrightarrow Y,
\qquad
(\Tx,\Ty,s,\phi,v)\longmapsto(\Tx+k\Ty,\Ty,s,\phi,a(s,\phi)v).
\]
Write $x=[d\Tx]$, $y=[d\Ty]$, and let $u$ be the degree-three generator from the last factor. Then
\[
\mathcal R=H^*(Y;\ZZ)
=\Lambda_\ZZ(x,y,\eta,u)\otimes\ZZ[\omega]/(\omega^2),
\]
where $x,y,\eta$ have degree one and $\omega,u$ have degrees two and three. The pullback is
\[
\Phi^*x=x+ky,\qquad \Phi^*u=u+m\omega\eta,
\qquad \Phi^*y=y,\quad \Phi^*\eta=\eta,\quad \Phi^*\omega=\omega.
\]

Use the product CW structure on $Y$, whose cellular cochain differential is zero, and choose a cellular approximation of $\Phi$ preserving the subcomplex over the basepoint of $S^2$. Its cellular cochain map is determined by $\Phi^*$, so the mapping torus has the integral cochain model
\[
\mathcal R\oplus d\vartheta\,\mathcal R,\qquad
D_R(a+d\vartheta\,b)=d\vartheta\,Ra,\qquad R=\Phi^*-1.
\]
Here $d\vartheta$ denotes the integral circle cochain. This model retains the filtration from the base $S^1\times S^2$.

Extend $\Ttw$ to the degree-zero derivation of $\mathcal R$ with $\Ttw x=y$, zero on the other generators, and define a degree-zero derivation $\mathsf U$ by $\mathsf Uu=\omega\eta$, zero on the other generators. They commute and satisfy $\Ttw^2=\mathsf U^2=0$. Hence, with $\mathcal N=k\Ttw+m\mathsf U$,
\begin{equation}\label{eq:integral-wang}
R=(1+k\Ttw)(1+m\mathsf U)-1
=k\Ttw+m\mathsf U+km\Ttw\mathsf U
=\mathcal N(1+m\mathsf U).
\end{equation}
The last factor is an integral automorphism with inverse $1-m\mathsf U$ and determinant one in every degree. Thus
\[
a+d\vartheta\,b\longmapsto(1+m\mathsf U)a+d\vartheta\,b
\]
is an integral chain isomorphism from $D_R$ to $D_{\mathcal N}$. Its correction raises the $S^2$ filtration degree by two, so it induces the identity on fibre cohomology.

Under $x\leftrightarrow\Tgen1$, $y\leftrightarrow\Tgenp$, $\eta\leftrightarrow u_1$ and $u\leftrightarrow u_3$, the differential $D_{\mathcal N}$ is precisely the finite model of \eqref{eq:monodromy-model} after cancelling the contractible pairs of the circle complexes in $V_m$. This cancellation is integral: choosing a spanning tree in the $m$-vertex circle eliminates $m-1$ pairs with coefficients $\pm1$ in each sphere degree. The remaining primitive classes are the $u_q$ above, and the higher map is $Q_*u_3=mu_1$. These comparisons preserve the base filtration and the fibrewise integral classes used to normalize $I$. By naturality of the spectral-sequence determinant maps, their determinant identification after realification agrees with $I_*$ up to sign. This supplies the integral model for the lattice calculations below.

For completeness, the Wang sequence is
\[
0\longrightarrow\operatorname{coker}R_{q-1}
\longrightarrow H^q(\Xkm;\ZZ)
\longrightarrow\ker R_q\longrightarrow0.
\]
Its right-hand term is free, so the sequence splits as a sequence of abelian groups. By \eqref{eq:integral-wang}, $R_q$ and $\mathcal N_q$ have the same Smith normal form. The only block mixing $k$ and $m$ is in degree four. In the basis $(xu,yu,x\omega\eta,y\omega\eta)$ it is
\[
\begin{pmatrix}
0&0&0&0\\
k&0&0&0\\
m&0&0&0\\
0&m&k&0
\end{pmatrix}.
\]
Its entries have greatest common divisor $\gkm$, and its nonzero $2\times2$ minors have greatest common divisor $\gkm^2$. Its nonzero Smith factors are therefore $\gkm,\gkm$. The remaining nonzero blocks are multiplication by $k$ in degrees $1,2,3,5,6,7$ and by $m$ in degrees $3,5$. The Wang sequence now gives the Betti numbers $1,3,5,6,7,7,6,5,3,1$ in degrees $0,\dots,9$, and
\[
\operatorname{tors}H^q(\Xkm;\ZZ)=
\begin{cases}
\ZZ/k,&q=2,3,7,8,\\
\ZZ/k\oplus\ZZ/m,&q=4,6,\\
(\ZZ/\gkm)^2,&q=5,\\
0,&\text{otherwise}.
\end{cases}
\]
On $\Lk\otimes\langle u_3,u_1\rangle$ the integral Leray lattice and integral cohomology modulo torsion differ. In degree four, $\Tgenp u_3$ lifts to the real cocycle $\Tgenp u_3-(m/k)\areaS\Tgen1u_1$, and the smallest integral multiple of this cocycle is $k/\gkm$ times it. So cohomology has index $k/\gkm$ in the Leray lattice. In degree five, the relation $\Akm(\Tgen1u_3)=k\,d\vartheta\,\Tgenp u_3+m\nuB\Tgen1u_1$ makes $\nuB\Tgen1u_1$, modulo torsion, equal to $k/\gkm$ times a primitive class. So the Leray lattice has index $k/\gkm$ in cohomology. Since odd degrees enter the determinant line inversely, both contribute the factor $k/\gkm$.

\emph{The determinant norm.}
Define $\ell_\ZZ(E_2)$ and $\ell_\ZZ(H)$ as in \S\ref{ssec:doublecover}. The flat bundle $H_\partial$ has a unimodular integral lattice, so by \cite{BZ} the Ray--Singer norm of $\ell_\ZZ(E_2)$ is the alternating product of the orders of the torsion of $H(B,H(Z;\ZZ))$. The groups $\ZZ/k$ lie in total degrees $2,4$, then $3,5$, then $5,7$, then $6,8$, for the rows $u_0,u_1,u_3,u_4$ of $\Lk$. Hence $\|\ell_\ZZ(E_2)\|_{\RS}=1$. The two rank-one differentials contribute $m^{-1}$ each, as in \S\ref{ssec:doublecover}, so with the lattice index above,
\[
\kappa\bigl(\ell_\ZZ(E_2)\bigr)=\pm\Bigl(\frac{\gkm}{km}\Bigr)^2\ell_\ZZ(H).
\]
Theorem~B therefore yields
\[
\|\ell_\ZZ(H)\|_{\RS(B,V,\Akm)}=\lkm^2.
\]
Independently, the torsion groups above give
\[
\prod_q|\operatorname{tors}H^q(\Xkm;\ZZ)|^{(-1)^q}
=k\cdot k^{-1}\cdot km\cdot\gkm^{-2}\cdot km\cdot k^{-1}\cdot k=\lkm^2,
\]
so classical Cheeger--M\"uller gives the same norm on $\Xkm$, as Corollary~\ref{cor:morse} requires. In Theorem~C the first sum vanishes, since each parity of $H_\partial$ consists of six rows with Betti numbers $1,1,1,1$ over $B$, and so $\chi'(\Akm)=-2$.

For $k=0$, $X_{0,m}=T^2\times\varpi^*X_m$, and $d_3$ has two further ranks, on the rows of $\Tgen1$ and $\Tgenp$, from even to odd total degree. The norm is then $1$ and $\chi'=0$, as the product formula with $\chi(T^2)=0$ requires. The Dehn twist removes those two ranks.

\section*{Acknowledgements}
The ideas in this paper were developed with assistance from Claude Opus 5.5 (Anthropic), which was also used for typesetting and drafting.


\begin{thebibliography}{99}

\bibitem{AC} C.~Arias Abad, M.~Crainic, \emph{Representations up to homotopy of Lie algebroids}, J.~Reine Angew.\ Math.\ \textbf{663} (2012), 91--126. Numbered references follow arXiv:0901.0319v2.

\bibitem{AS} C.~Arias Abad, F.~Sch\"atz, \emph{Reidemeister torsion for flat superconnections}, J.~Homotopy Relat.\ Struct.\ \textbf{9} (2014), 579--606. Numbered references follow arXiv:1108.5103v1.

\bibitem{AS2} C.~Arias Abad, F.~Sch\"atz, \emph{The $A_\infty$ de Rham theorem and integration of representations up to homotopy}, Int.\ Math.\ Res.\ Not.\ IMRN \textbf{2013} (16), 3790--3855. Numbered references follow arXiv:1011.4693v2.

\bibitem{BGV} N.~Berline, E.~Getzler, M.~Vergne, \emph{Heat kernels and Dirac operators}, Springer, 1992.

\bibitem{BL} J.-M.~Bismut, J.~Lott, \emph{Flat vector bundles, direct images and higher real analytic torsion}, J.~Amer.\ Math.\ Soc.\ \textbf{8} (1995), 291--363.

\bibitem{BZ} J.-M.~Bismut, W.~Zhang, \emph{An extension of a theorem by Cheeger and M\"uller}, Ast\'erisque \textbf{205} (1992).

\bibitem{BS} J.~Block, A.~M.~Smith, \emph{The higher Riemann--Hilbert correspondence}, Adv.\ Math.\ \textbf{252} (2014), 382--405.

\bibitem{Ch} J.~Cheeger, \emph{Analytic torsion and the heat equation}, Ann.\ of Math.\ \textbf{109} (1979), 259--322.

\bibitem{CHL} J.~Chuang, J.~Holstein, A.~Lazarev, \emph{Maurer--Cartan moduli and theorems of Riemann--Hilbert type}, Appl.\ Categ.\ Structures \textbf{29} (2021), 685--728. Numbered references follow arXiv:1802.02549v3.

\bibitem{CL} J.~Chuang, A.~Lazarev, \emph{On the perturbation algebra}, J.~Algebra \textbf{519} (2019), 130--148. Numbered references follow arXiv:1703.05296v2.

\bibitem{Cr} M.~Crainic, \emph{On the perturbation lemma, and deformations}, arXiv:math/0403266.

\bibitem{Fa95} M.~Farber, \emph{Singularities of the analytic torsion}, J.~Differential Geom.\ \textbf{41} (1995), 528--572.

\bibitem{Fa} M.~Farber, \emph{Combinatorial invariants computing the Ray--Singer analytic torsion}, Differential Geom.\ Appl.\ \textbf{6} (1996), 351--366.

\bibitem{Fo} R.~Forman, \emph{Hodge theory and spectral sequences}, Topology \textbf{33} (1994), 591--611.

\bibitem{Fd} D.~S.~Freed, \emph{Reidemeister torsion, spectral sequences, and Brieskorn spheres}, J.~Reine Angew.\ Math.\ \textbf{429} (1992), 75--90.

\bibitem{Gi} P.~B.~Gilkey, \emph{Invariance theory, the heat equation, and the Atiyah--Singer index theorem}, Publish or Perish, 1984.

\bibitem{Go1} S.~Goette, \emph{Morse theory and higher torsion invariants I}, arXiv:math/0111222.

\bibitem{Go2} S.~Goette, \emph{Morse theory and higher torsion invariants II}, arXiv:math/0305287.

\bibitem{Go3} S.~Goette, \emph{Torsion invariants for families}, Ast\'erisque \textbf{328} (2009), 161--206; arXiv:0804.3020.

\bibitem{GSM} A.~Gracia-Saz, R.~A.~Mehta, \emph{Lie algebroid structures on double vector bundles and representation theory of Lie algebroids}, Adv.\ Math.\ \textbf{223} (2010), 1236--1275; arXiv:0810.0066.

\bibitem{Ka} T.~Kato, \emph{Perturbation theory for linear operators}, Springer, 1995.

\bibitem{KM} F.~F.~Knudsen, D.~Mumford, \emph{The projectivity of the moduli space of stable curves I}, Math.\ Scand.\ \textbf{39} (1976), 19--55.

\bibitem{LST} W.~L\"uck, T.~Schick, T.~Thielmann, \emph{Torsion and fibrations}, J.~Reine Angew.\ Math.\ \textbf{498} (1998), 1--33.

\bibitem{Ma} X.~Ma, \emph{Functoriality of real analytic torsion forms}, Israel J.\ Math.\ \textbf{131} (2002), 1--50.

\bibitem{MW1} V.~Mathai, S.~Wu, \emph{Analytic torsion for twisted de Rham complexes}, J.~Differential Geom.\ \textbf{88} (2011), 297--332; arXiv:0810.4204.

\bibitem{MW2} V.~Mathai, S.~Wu, \emph{Analytic torsion of $\ZZ_2$-graded elliptic complexes}, in: Noncommutative geometry and global analysis, Contemp.\ Math.\ \textbf{546}, Amer.\ Math.\ Soc., 2011, 199--212; arXiv:1001.3212.

\bibitem{Mi} R.~Mickler, \emph{Twisted analytic torsion and adiabatic limits}, Proc.\ Amer.\ Math.\ Soc.\ \textbf{143} (2015), no.~12, 5455--5469. Numbered references follow arXiv:1311.6788v1.

\bibitem{Mu} W.~M\"uller, \emph{Analytic torsion and R-torsion of Riemannian manifolds}, Adv.\ Math.\ \textbf{28} (1978), 233--305.

\bibitem{Q} D.~Quillen, \emph{Superconnections and the Chern character}, Topology \textbf{24} (1985), 89--95.

\bibitem{RS} D.~B.~Ray, I.~M.~Singer, \emph{R-torsion and the Laplacian on Riemannian manifolds}, Adv.\ Math.\ \textbf{7} (1971), 145--210.

\bibitem{So} O.~N.~Sorto, \emph{Twisted Milnor torsion for finite group actions}, PhD thesis, University of Auckland, 2022.

\end{thebibliography}
\end{document}